\documentclass{article}
\usepackage[english]{}
\usepackage[utf8]{inputenc}

\usepackage{amsmath}
\usepackage{amssymb}
\usepackage{graphicx}
\usepackage{color}
\usepackage{amsfonts}
\usepackage{setspace}
\usepackage[absolute,overlay]{textpos}
\usepackage{tikz}
\usepackage{multirow}
\usepackage{caption}
\usepackage{subcaption}
\usepackage{authblk}
\usepackage{placeins}

\usepackage[right]{lineno}
\usepackage{hyperref}
\usepackage{comment}
\title{Dynamic nonlinear multicontinuum homogenization of systems with intrinsically evolving microstructure}
\date{}

\author[a]{Mohammed Al Kobaisi}
\author[a]{Dmitry Ammosov}
\author[b]{Yalchin Efendiev}
\author[c]{Wing Tat Leung}
\author[b]{Buzheng Shan}

\affil[a]{Chemical and Petroleum Engineering Department, Khalifa University of Science and Technology, Abu Dhabi, 127788, UAE}
\affil[b]{Department of Mathematics, Texas A\&M University, College Station, TX 77843, USA}
\affil[c]{Department of Mathematics, City University of Hong Kong, Hong Kong}

\begin{document}
\maketitle

\renewcommand{\thefootnote}{}%
\footnotetext{E-mail addresses: mohammed.alkobaisi@ku.ac.ae (Mohammed Al Kobaisi), dmitrii.ammosov@ku.ac.ae (Dmitry Ammosov), efendiev@math.tamu.edu (Yalchin Efendiev), wtleung27@cityu.edu.hk (Wing Tat Leung), shanbzh@tamu.edu (Buzheng Shan).}
\addtocounter{footnote}{-1}%

\begin{abstract}
In this paper, we propose a multicontinuum homogenization approach for nonlinear problems involving dynamically evolving multiscale media. The main idea of the proposed approach is that one of the fine-scale variables defines continua. It allows us to formulate macroscopic variables and derive new macroscopic models for nonlinear problems, where coefficients can depend on fine-scale functions. As an example, we consider a fingering problem and employ the fine-scale concentration field to define continua. We consider both Galerkin and mixed multicontinuum modeling approaches. In the former, the multicontinuum theory is applied to the pressure and concentration fields; in the latter, it is also applied to the velocity field. In both approaches, we provide multicontinuum expansions, formulate cell problems, and derive the corresponding macroscopic models. We present numerical results for model problems of gravity-driven fingering, viscous fingering, and interface flattening driven by high-contrast flow. The results show that the macroscopic models, derived with the proposed approach, can provide an accurate representation of the coarse-scale solutions.

\end{abstract}

\section{Introduction}

Many practical problems have multiscale nature. Multiple scales can
have complex spatial and temporal configurations. These complex features
can exhibit themselves as highly anisotropic features, 
which are highly connected, e.g., viscous fingering \cite{homsy1987viscous, shokri2018saffman}. These multiscale
features evolve in time, and their modeling on a coarse grid is challenging.

Many advanced multiscale methods have been proposed to capture complex
heterogeneities for various types of processes. As an example,
we focus on CEM-GMsFEM \cite{chung2018constraint} and also mention several
other multiscale methods that can alternatively be used to solve
challenging multiscale problems \cite{chen2003mixed, hou1997multiscale, jenny2003multi, efendiev2009multiscale, efendiev2013generalized, maalqvist2014localization,chung2015mixed, chaabi2024algorithmic, peter2009multiscale, henning2013oversampling, owhadi2014polyharmonic,fu2025wavelet, xie2026robust, panasenko2018multicontinuum, roberts2014dynamical}. For simplicity
in this section, we consider 
a typical equation
\[
u_t + L(u)=f,
\]
where $u$ is the solution, and $L(u)$ is a differential operator, 
e.g., $L(u)=-div(\kappa(x)\nabla u)$, with $\kappa(x)$ being a
heterogeneous field. 
In CEM-GMsFEM, the solution is sought in terms of multiscale basis functions
\[
u^{\text{CEM}} = \sum_{i,j} u_i^j(t) \widetilde{\phi}_i^j(x),
\]
where $\widetilde{\phi}_i^j$ is a multiscale basis function at the coarse node $j$
($i$ being the numbering of basis functions). These basis functions, supported
in local domains, are
computed using local spectral problems and energy minimization.  It was shown that $u^{\text{CEM}}$ has a first order
convergence.

In recent works, we propose multicontinuum homogenization methods \cite{efendiev2023multicontinuum, chung2024multicontinuum, leung2024some}.
These methods employ multiscale basis functions such that the
coefficients $u_i^j(t)$ are continuous with respect to space
variables, which allows formulating macroscopic PDEs, i.e.,
$u_i^j(t)\approx u_i(x,t)$, where the index $j$ is replaced by
a spatial variable. We seek the approximation of the solution in the form
\[
u^{\text{MH}} = \sum_{i} \phi_i(x) u_i(x, t) +
\sum_{i,m} \phi_i^{m}(x) \nabla_m u_i(x, t),
\]
where $\nabla_m$ refers to $\frac{\partial}{\partial x_m}$, $\phi_i(x)$ are multicontinuum basis functions
for averages, and $\phi_i^{m}(x)$ are for gradients. 
The latter is a gradient part of the basis functions $\widetilde{\phi}_i^j$.
These
basis functions solve local problems
and  use local multiscale basis functions,
denoted by $\psi_j(x)$,
 as constraints
\[
\int \phi_i(x)\psi_k(x)=\delta_{ik},\ \ \int \phi_i^{m}(x)\psi_k(x)=\delta_{ik} (x_m-x_m^0), 
\]

Next, we briefly motivate multicontinuum homogenization.
We denote  $U$
to consist of all macroscopic variables, and
$\Pi^{\text{MS}}$ is a projection of macroscopic variables; then
\begin{equation*}
    \Pi^{\text{MS}}(u)\approx\Pi^{\text{MS}}(\psi_iU_i).
\end{equation*}
If we consider $U_i$'s are smooth and $\Pi^{\text{MS}}$ to be local (in a certain sense), then we further obtain
\begin{equation}
\begin{split}
\label{eq:motiv}
    \Pi^{\text{MS}}(u) 
    &\approx \Pi^{\text{MS}}(\psi_i{\bf 1}U_i(x^*_0)+\psi_i(x-x_0^*)\cdot\nabla U_i (x^0))\\
    &\approx \Pi^{\text{MS}}(\psi_i {\bf 1})U_i(x^*_0)+\Pi^{\text{MS}}(\psi_i(x-x_0^*))\cdot\nabla U_i(x^*_0)
    =\phi_iU_i(x^*_0) + \phi_i^m\nabla_m U_i(x^*_0),
\end{split}
\end{equation}
in the region near $x^*_0$ where $\phi_i := \Pi^{\text{MS}}(\psi_i{\bf 1})$, and $\phi_i^m := (\Pi^{\text{MS}}(\psi_i(x-x_0^*)))_m$. Note that from here on, we assume summation over repeated indices.

\begin{figure}[hbt!]
\centering
\includegraphics[width=0.49\textwidth]{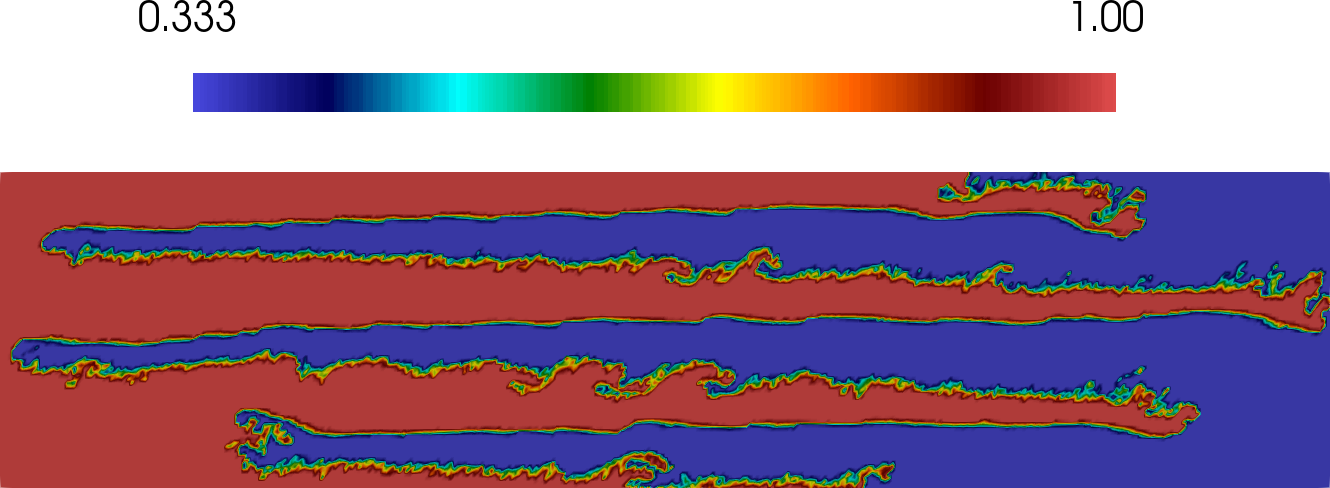}
\caption{A snapshot of the solution}
\label{fig:snapshots}
\end{figure}

One of our goals with multicontinuum homogenization methods is to 
derive macroscopic equations for complex nonlinear problems. In these
problems, macroscopic properties are nonlinear functions. We consider
one example,
\begin{equation}
\label{eq:general}
\begin{split}
-div(\lambda(c,x) (\nabla p-c e_1)) =f,\\
c_t + \nabla \cdot(v c) =0,\ v = -\lambda(c,x) (\nabla p-c e_1),
\end{split}
\end{equation}
where $p$ is the pressure field, and $c$ is the concentration field.
Here, for simplicity, we choose $\lambda=1$, which is a well-known
gravity-driven fingering problem. This problem is known to be ill-posed and 
strongly depends on regularization \cite{petrova2025propagating}. In our paper, for computing the concentration 
at the $(n + 1)$th time step, we use the velocity from the $n$th time step, i.e., based on $\lambda(c^{n}, x)$.
We present a time snapshot of the solution in Figure \ref{fig:snapshots}.

The main idea of the proposed approach is that one of the variables,
 $c(x,t)$, 
defines continua. In particular, we define each continua as
\[
I_k=\{ c_k\leq c < c_{k+1} \}.
\]
For each $I_k$, we define a characteristic function $\psi_k$, 
\[
\psi_k(c)=\{1 \text{ if } c\in I_k;0 \text{ otherwise } \}.
\]
This
identifies the continuum at each time. 
In previous approaches, the continua are
stationary and 
defined by high and low contrast values of the permeability,
which can also be added here. 

The main objective is to propose a general formalism for
macroscopic modeling. 
Given each continuum, we define an expansion of the solution         
\begin{equation}
\begin{split}
p=\phi^p_iP_i + \phi_i^{p,m} \nabla_m P_i + \phi^{pc}_iC_i +\phi^{pc,m}_i\nabla_m C_i\\
c=\phi^{c}_iC_i+\phi^{c,m}_i\nabla_m C_i. \\
\end{split}
\end{equation}
Here, we linearize the process around the concentration. As we mentioned earlier, this expansion is motivated by multiscale basis functions (see (\ref{eq:motiv})). In the above problems, $\phi^{pc,m}_i$ and $\phi^{c,m}_i$ are small, and we neglect them in our derivations and numerical simulations. The basis functions in multicontinuum models
are constructed using the constraints as formulated below. 
These constraints are chosen such that the macroscopic variables
are smooth functions.

{\bf Cell problems.} The cell solutions $\phi^p_i$ and $\phi^{p,m}_i$
follow the standard elliptic equation construction. For 
$\phi^{pc}_i$ and $\phi^{c}_i$, we have the cell
problems formulated in (\ref{eq:constraint2}). 
The equations for 
$\phi^{pc,m}_i$ and $\phi^{c,m}_i$ are the same as (\ref{eq:constraint2})
with constraints on gradients 
$\int_{R_\omega^{l}} \phi^{c,m}_i\psi_j^l = \delta_{ij}\int_{R_\omega^{l}} (x_m-x_m^0)\psi_j^l.$
Using this expansion, we arrive at the following macroscopic equations
(or their modifications)
\begin{equation}
\begin{split}
-\nabla_n(\alpha_{ij}^{mn} \nabla_m P_j)+\beta_{ij} P_j +\nabla_m (\gamma^m_{ij}C_j) + \gamma_{ij} C_j  =f_i^p\\
(C_i)_t + \zeta_{ij} C_j + \eta_{ij}^m\nabla_m C_j =f_i^c,
\end{split}
\end{equation}
where the macroscopic coefficients depend on concentrations and will be defined
later. We note that depending on multicontinuum expansion, we can 
arrive at different formulations of macroscopic equations that involve
more terms; however, many coefficients in the resulting equations can become
negligle (cf. \cite{ammosov2025multicontinuum_poro}).

In our approach, we consider both Galerkin and mixed multicontinuum
modeling. The difference is that in the former, we apply
multicontinuum theory for pressure only, while in the latter,
we apply the multicontinuum theory for pressure and velocity.

Our main contributions are the following.

\begin{itemize}

\item
Our approach allows formulating 
macroscopic variables and deriving {\it new}
macroscopic models for nonlinear
problems, where coefficients can depend on fine-grid functions. 

\item  
One can attempt to approximate the macroscopic coefficients using coarse-grid
quantities or compute based on local solutions using machine learning
techniques

\item
Our main objective is to show that the derived macroscopic models 
can provide
an accurate representation of the coarse-grid solution.

\end{itemize}

We present numerical results for model problems in dual- and triple-continuum media. Specifically, we consider cases where the continua are determined by gravitational effects and high-contrast properties. Among the former, we examine gravity-driven fingering problems, while for the latter, we study problems of viscous fingering and interface flattening driven by high-contrast flow. Note that we apply the Lagrangian-Eulerian particle-mesh method based on the LeoPart library \cite{maljaars2021leopart} to avoid numerical diffusion in gravity-driven fingering problems on a fine grid. Depending on the problem type, we apply either mixed or Galerkin multicontinuum modeling approaches on a coarse grid. Overall, the results demonstrate that our proposed multicontinuum models can approximate the reference average solutions in dynamic multicontinuum media with high accuracy.

The paper has the following structure. In Section \ref{sec:proposed_method}, we present our proposed dynamic multicontinuum modeling approach. We describe both Galerkin and mixed multicontinuum approaches and derive the corresponding macroscopic models. Section \ref{sec:numerical_results} presents numerical results for model problems in dual- and triple-continuum media. We consider both gravity-determined and contrast-determined continua cases. Finally, we present conclusions in Section \ref{sec:conclusions}.

\section{Proposed method}\label{sec:proposed_method}

In this section, we present the details of our approach for 
problems (\ref{eq:general}).  
We will consider two different derivations for pressure
equations, with Galerkin and mixed
formulations. Both approaches have their advantages in simulations.

\subsection{Coarse and fine grids}

\begin{figure}[hbt!]
\centering
\includegraphics[width=0.5\textwidth]{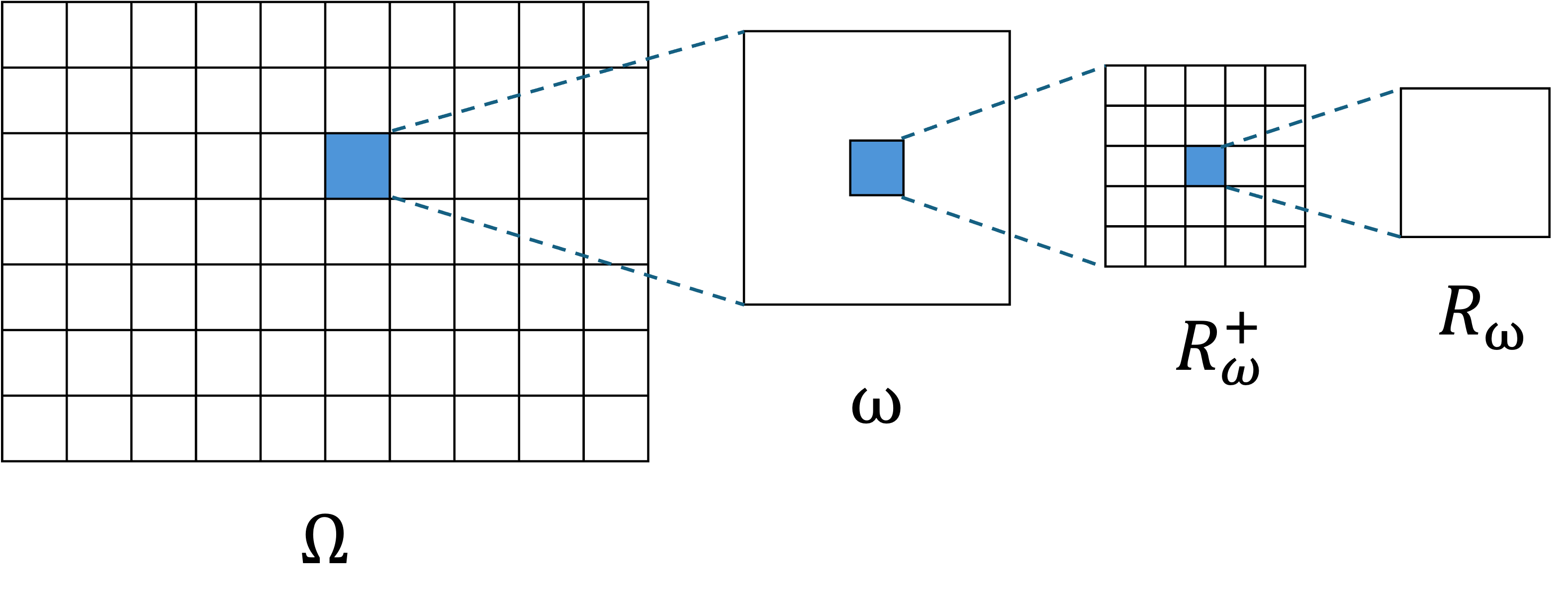}
\caption{Illustration of the domain $\Omega$, coarse block $\omega$, oversampled RVE $R_\omega^+$, and RVE $R_\omega$}
\label{fig:mh_scheme}
\end{figure}

Let $\Omega$ be a computational domain, partitioned into coarse blocks $\omega$. Inside each coarse block, we assume that there is a Representative Volume Element (RVE) $R_\omega$ that can represent the whole $\omega$. In addition, we suppose that we can construct an oversampled RVE $R_\omega^+$ as an extension of $R_\omega$ with layers of neighbour RVEs, i.e., $R_\omega^+ = \cup_l R_\omega^l$, where $R_\omega^{l_0} = R_\omega$ (see Figure \ref{fig:mh_scheme}). Finally, we assume that we can identify $N$ components (continua) in each RVE that depend on $c$ as described earlier. We introduce a characteristic function $\psi_j(c)$ for each continuum $j$, which equals $1$ in the continuum $j$ and $0$ otherwise.

Based on the characteristic functions, we define macroscopic variables, which represent the average solutions in the corresponding subregions
\begin{equation}\label{eq:macroscopic_variable_definitions}
P_{i}(x_\omega) = \frac{\int_{R_\omega} p \psi_i}{\int_{R_\omega} \psi_i}, \quad
C_i(x_\omega) = \int_{R_\omega} c \psi_i,
\end{equation}
where $x_\omega$ is a point in $R_\omega$. Note that we assume that $C_{i}$ and $P_{i}$ are smooth functions. In general, we can have different number of continua for pressure and concentration.

\subsection{Galerkin approach}\label{sec:galerkin_approach}

In the present approach, we consider the following expansion
\begin{equation}
\begin{split}
p=\phi^p_iP_i + \phi_i^{p,m} \nabla_m P_i + \phi^{pc}_iC_i +\phi^{pc,m}_i\nabla_m C_i\\
c=\phi^{c}_iC_i+\phi^{c,m}_i\nabla_m C_i. \\
\end{split}
\end{equation}

The cell solutions $\phi^p_i$ and $\phi^{p,m}_i$
follow the standard elliptic equation construction based on the following cell problems.
\begin{equation}
\label{eq:cell_grad_aver}
\begin{split}
\int_{R_\omega^+}\lambda \nabla \phi^{p, m}_i\cdot \nabla q -  \sum_{j,l} {\beta_{ij}^{ml}\over \int_{R_\omega^l}\psi_j^l} \int_{R_\omega^l}\psi_j^l  q  =0\\
\int_{R_\omega^l}  \phi^{p, m}_i \psi_j^l = \delta_{ij} \int_{R_\omega^l} (x_m-\tilde{x}_{mj})\psi_j^l ,\\
\int_{R_\omega^{l_0}} (x_m-\tilde{x}_{mj})\psi_j^{l_0} =0\ \text{condition that defines $\tilde{x}$}, \\
\end{split}
\end{equation}
where $m$ refers to the coordinate direction,
and
\begin{equation}
\label{eq:cell_aver_aver}
\begin{split}
\int_{R_\omega^+}\lambda \nabla \phi_i^p\cdot \nabla q -
\sum_{j,l} {\beta_{ij}^l\over \int_{R_\omega^l}\psi_j^l} \int_{R_\omega^l}\psi_j^l  q =0\\
\int_{R_\omega^l}  \phi_i^p \psi_j^l = \delta_{ij} \int_{R_\omega^l} \psi_j^l.
\end{split}
\end{equation}
 For 
$\phi^{pc}_i$ and $\phi^{c}_i$, we have the following cell problems
\begin{equation}
\label{eq:constraint2}
\begin{split}
-div(\lambda(\cdot,\cdot) (\nabla\phi^{pc}_i -\phi^{c}_i e_1)) =0\\
\int_{R_\omega^{+}} v \phi^{c}_i\nabla \xi - \sum_{j,l}{\beta_{ij}^l \over \int_{R_\omega^{l}} \psi_j^l}\int_{R_\omega^{l}} \xi \psi_j^l =0\\
\int_{R_\omega^{l}} \phi^{c}_i\psi_j^l = \delta_{ij}\int_{R_\omega^{l}} \psi_j^l.
\end{split}
\end{equation}

Note that we will neglect $\phi^{pc,m}_i$ and $\phi^{c,m}_i$ in our derivations and numerical simulations. 
However, if we use $\phi^{pc,m}_i$ and $\phi^{c,m}_i$, then the equations for them are the same as (\ref{eq:constraint2})
with constraints on gradients $\int_{R_\omega^{l}} \phi^{c,m}_i\psi_j^l = \delta_{ij}\int_{R_\omega^{l}} (x_m-\tilde{x}_m)\psi_j$.

Using the multicontinuum expansion, we arrive at the following macroscopic equations
(see the derivation in Section \ref{sec:derG})
\begin{equation}
\label{eq:macG}
\begin{split}
-\nabla_n (\alpha_{ij}^{mn} \nabla_m P_j)+\beta_{ij} P_j - \nabla_m (\gamma^m_{ij}C_j) + \gamma_{ij} C_j  =f_i^p\\
(C_i)_t + \zeta_{ij} C_j + \eta_{ij}^m\nabla_m C_j =f_i^c,
\end{split}
\end{equation}
where the macroscopic coefficients depend on concentrations.

In our numerical
simulations, we take $\phi_i^{c}$ to be piecewise constant functions 
based on values of the concentration as discussed above. 
We consider
\[
c=\phi_i^{c}(x) C_i(x,t),
\]
where $\phi_i^{c}(x) =\chi_i \psi_i$, $\chi_i = (\int_{R_\omega} \psi_i)^{-1}$.

\subsubsection{Derivation of the macroscopic model}\label{sec:derG}
We first present the derivation of the pressure equations and then
the derivation of the concentration equations.

\paragraph{Pressure equations}

For simplicity, we ignore the terms with multicontinuum concentration gradients, since they are small. Therefore, we consider the following multicontinuum expansions
\begin{equation}\label{eq:mc_expansions_simpl}
\begin{split}
p \approx \phi^p_iP_i + \phi_i^{p,m} \nabla_m P_i + \phi^{pc}_i C_i,\\
c \approx \phi_i^c C_i.
\end{split}
\end{equation}

Next, we suppose homogeneous Dirichlet boundary conditions and obtain the following variational formulation of the pressure equation
\begin{equation}\label{eq:var_form_p}
\int_\Omega \lambda(c,x) \nabla p \cdot \nabla q - \int_\Omega \lambda(c,x) c e_1 \cdot \nabla q = \int_\Omega f q, \quad \forall q \in H_0^1(\Omega).\\
\end{equation}

Then, we use the RVE property to approximate the whole domain and substitute the expansion \eqref{eq:mc_expansions_simpl} into this variational formulation
\begin{equation}\label{eq:mc_pressure_RVE_variational}
\begin{split}
\int_\Omega \lambda(c,x) \nabla p \cdot \nabla q 
- \int_\Omega \lambda(c,x) c e_1 \cdot \nabla q \approx \\
\sum_\omega \frac{|\omega|}{|R_\omega|} \left\{ \int_{R_\omega} \lambda(c,x) \nabla p \cdot \nabla q - 
\int_{R_\omega} \lambda(c,x) c e_1 \cdot \nabla q \right\} \approx\\
\sum_\omega \frac{|\omega|}{|R_\omega|} \left\{ \int_{R_\omega} \lambda(c,x) \nabla (\phi^p_iP_i + \phi_i^{p,m} \nabla_m P_i + \phi^{pc}_i C_i) \cdot \nabla q - 
\int_{R_\omega} \lambda(c,x) \phi_i^c C_i e_1 \cdot \nabla q \right\} =\\
\sum_\omega \frac{|\omega|}{|R_\omega|} f q.
\end{split}
\end{equation}

Let us consider the first integral term. To obtain the multicontinuum pressure equations, we consider only the pressure parts of the test function's expansion, i.e., we take $q \approx \phi^p_j Q_j + \phi_j^{p,n} \nabla_n Q_j$. Moreover, we employ our assumption about the smoothness of macroscopic variables $P_i$ and $C_i$ and assume that their variations are minor compared to $\phi_i^p$, $\phi_i^{p,m}$, and $\phi_i^{pc}$
\begin{equation}\label{eq:mc_pressure_variational_1}
\begin{split}
\int_{R_\omega} \lambda(c,x) \nabla (\phi^p_i P_i + \phi_i^{p,m} \nabla_m P_i + \phi^{pc}_i C_i) \cdot \nabla q \approx\\
\int_{R_\omega} \lambda(c,x) \nabla (\phi^p_i P_i + \phi_i^{p,m} \nabla_m P_i + \phi^{pc}_i C_i) \cdot \nabla (\phi^p_j Q_j + \phi_j^{p,n} \nabla_n Q_j) \approx \\
P_i Q_j \int_{R_\omega} \lambda(c,x) \nabla \phi^p_i \cdot \nabla \phi^p_j +
\nabla_m P_i Q_j \int_{R_\omega} \lambda(c,x) \nabla \phi^{p,m}_i \cdot \nabla \phi^p_j +\\
C_i Q_j \int_{R_\omega} \lambda(c,x) \nabla \phi^{pc}_i \cdot \nabla \phi^p_j +
P_i \nabla_n Q_j \int_{R_\omega} \lambda(c,x) \nabla \phi^p_i \cdot \nabla \phi^{p, n}_j +\\
\nabla_m P_i \nabla_n Q_j \int_{R_\omega} \lambda(c,x) \nabla \phi^{p,m}_i \cdot \nabla \phi^{p, n}_j +
C_i \nabla_n Q_j \int_{R_\omega} \lambda(c,x) \nabla \phi^{pc}_i \cdot \nabla \phi^{p, n}_j = \\
|R_\omega| (
\beta_{ij} P_i Q_j +
\alpha_{ij}^m \nabla_m P_i Q_j + 
\overline{\gamma}_{ij} C_i Q_j +
\overline{\alpha}_{ij}^n P_i \nabla_n Q_j + 
\alpha_{ij}^{nm} \nabla_m P_i \nabla_n Q_j +
\overline{\gamma}_{ij}^n C_i \nabla_n Q_j ),
\end{split}
\end{equation}
where the effective properties are defined as follows
\begin{equation*}
\begin{gathered}
\beta_{ij} = \frac{1}{|R_\omega|} \int_{R_\omega} \lambda(c,x) \nabla \phi^p_i \cdot \nabla \phi^p_j, \quad
\alpha_{ij}^{nm} = \frac{1}{|R_\omega|} \int_{R_\omega} \lambda(c,x) \nabla \phi^{p,m}_i \cdot \nabla \phi^{p, n}_j, \\
\alpha_{ij}^m = \frac{1}{|R_\omega|} \int_{R_\omega} \lambda(c,x) \nabla \phi^{p,m}_i \cdot \nabla \phi^p_j, \quad
\overline{\alpha}_{ij}^n = \frac{1}{|R_\omega|} \int_{R_\omega} \lambda(c,x) \nabla \phi^p_i \cdot \nabla \phi^{p, n}_j, \\
\overline{\gamma}_{ij} = \frac{1}{|R_\omega|} \int_{R_\omega} \lambda(c,x) \nabla \phi^{pc}_i \cdot \nabla \phi^p_j, \quad
\overline{\gamma}_{ij}^n = \frac{1}{|R_\omega|} \int_{R_\omega} \lambda(c,x) \nabla \phi^{pc}_i \cdot \nabla \phi^{p, n}_j.
\end{gathered}
\end{equation*}

The second integral term in \eqref{eq:mc_pressure_RVE_variational} is approximated in a similar way
\begin{equation}\label{eq:mc_pressure_variational_2}
\begin{split}
-\int_{R_\omega} \lambda(c,x) \phi_i^c C_i e_1 \cdot \nabla q \approx
-\int_{R_\omega} \lambda(c,x) \phi_i^c C_i e_1 \cdot \nabla (\phi^p_j Q_j + \phi_j^{p,n} \nabla_n Q_j) =\\
-\int_{R_\omega} \lambda(c,x) \phi_i^c C_i e_1 \cdot \nabla \phi^p_j Q_j
-\int_{R_\omega} \lambda(c,x) \phi_i^c C_i e_1 \cdot \nabla \phi_j^{p,n} \nabla_n Q_j \approx\\
-C_i Q_j \int_{R_\omega} \lambda(c,x) \phi_i^c e_1 \cdot \nabla \phi^p_j
-C_i \nabla_n Q_j \int_{R_\omega} \lambda(c,x) \phi_i^c e_1 \cdot \nabla \phi_j^{p,n} =\\
-|R_\omega| ( \tilde{\gamma}_{ij} C_i Q_j +
\tilde{\gamma}_{ij}^n C_i \nabla_n Q_j ),
\end{split}
\end{equation}
where the effective properties are following
\begin{equation*}
\tilde{\gamma}_{ij} = \frac{1}{|R_\omega|} \int_{R_\omega} \lambda(c,x) \phi_i^c e_1 \cdot \nabla \phi^p_j, \quad
\tilde{\gamma}_{ij}^n = \frac{1}{|R_\omega|} \int_{R_\omega} \lambda(c,x) \phi_i^c e_1 \cdot \nabla \phi_j^{p,n}.
\end{equation*}

We sum the terms from \eqref{eq:mc_pressure_variational_1} and \eqref{eq:mc_pressure_variational_2} over blocks. 
Moreover, we neglect the second and fourth terms in \eqref{eq:mc_pressure_variational_1}, since their sum is small \cite{efendiev2023multicontinuum}. In this way, we obtain the following multicontinuum pressure equations
\begin{equation}
-\nabla_n (\alpha_{ij}^{mn} \nabla_m P_j) + \beta_{ij} P_j - \nabla_m (\gamma^m_{ij} C_j) + \gamma_{ij} C_j = f_i^p,
\end{equation}
where $\gamma^m_{ij} = \overline{\gamma}_{ji}^m - \tilde{\gamma}_{ji}^m$, $\gamma_{ij} = \overline{\gamma}_{ji} - \tilde{\gamma}_{ji}$, and $f_i^p = \frac{1}{|R_\omega|} \int_{R_\omega} f \phi^p_i$.

\paragraph{Concentration equations.}

We present a formal derivation of macroscopic equations, where
the concentration equation is multiplied by $\psi_k$.
Multiplying the transport equation by 
$\psi_k(c)$, we obtain
\begin{equation}\label{eq:multicontinuum_concentration_model}
\begin{split}
0=\int (c)_t\psi_k   + \int v \nabla c \psi_k =
\int (c\psi_k)_t -  \int c(\psi_k)_t   + \int_{\partial_{\psi_k^{ext}}} v\cdot n c -
 \int v  c \nabla \psi_k=\\
(C_k)_t  + \int_{\partial_{\psi_k^{ext}}} v\cdot n c - \int c(-v \nabla \psi_k) -
 \int v  c \nabla \psi_k=\\
(C_k)_t  + \int_{\partial R_\omega} v\cdot n c\psi_k.
\end{split}
\end{equation}
Next, we use the following expansions
\begin{equation}
c\approx  \phi_j^c C_j,\ \ 
v\approx \phi_j^v V_j + \phi_j^{vp} P_j - \phi_j^c C_j e_1
\end{equation}
to compute 
\begin{equation}
\begin{split}
\int_{\partial R_\omega} (\phi_j^v V_j + \phi_j^{vp} P_j - \phi_j^c C_j e_1) \cdot n \phi_m^c C_m \psi_k \approx \\
\int_{\partial R_\omega}  \left\{ \phi_j^v \phi_k^c (V_j\cdot n) C_k + (\phi_j^{vp} \cdot n \phi_k^c) P_j C_k + (e_1 \cdot n) \phi_k^c \phi_k^c C_k C_k \right\} =\\
\int_{R_\omega} \left\{ div(\alpha_{jk} V_j C_k) +\nabla_n (\beta_{jk}^n P_j C_k) + \nabla_x (\gamma_k C_k C_k)\right\},
\end{split}
\end{equation}
where $\alpha_{jk}$, $\beta_{jk}^n$, and $\gamma_{k}$ are appropriately defined. If we have one continua for pressure, then the second term disappears. Note that, here, we have used $\psi_i \phi_i^c = \psi_i \chi_i \psi_i = \chi_i \psi_i = \phi_i^c$.

We present another derivation in the Appendix \ref{sec:appendix}. Next, we discuss macroscopic equations for pressure field within the mixed framework.

\subsection{Mixed framework}

In this subsection, we consider the mixed multicontinuum framework. Since the derivation of multicontinuum concentration equations is the same as in Section \ref{sec:derG}, we derive only multicontinuum flow equations. The mixed formulation of the fine-scale problem \eqref{eq:general} is following
\begin{equation}
\label{eq:general1}
\begin{split}
\lambda(c,x)^{-1} v + \nabla p-c e_1=0,\\
div(v)=f,\\
c_t + \nabla \cdot(v c) =0.
\end{split}
\end{equation}

Note that we can represent this mixed flow model in the following form
\begin{equation}
A Y = B,
\end{equation}
where the operator $A$, the solution field vector $Y$, and the right-hand side vector $B$ are defined as follows
\begin{equation}
A = \begin{bmatrix}
\lambda^{-1} & \nabla & -e_1\\
div & 0 & 0
\end{bmatrix}, \quad
Y = \begin{bmatrix}
v\\
p\\
c
\end{bmatrix}, \quad
B = \begin{bmatrix}
0\\
f
\end{bmatrix}.
\end{equation}

Next, we introduce the following multicontinuum expansions of the solution fields
\begin{equation}\label{eq:mc_mixed_expansion_2}
\begin{split}
p=\phi_i^{pv} V_i + \phi_i^p P_i + \phi_i^{p,m} \nabla_m P_i + \phi_i^{pc} C_i,\\
v=\phi_i^v V_i + \phi_i^{vp} P_i + \phi_i^{vp,m} \nabla_m P_i + \phi_i^{vc} C_i,\\
c=\phi_i^c C_i.
\end{split}
\end{equation}
We found in our numerical simulations that the $\phi_i^{pc}$ terms are negligible in the pressure expansion.

We can obtain the multicontinuum basis functions $\phi_i^v = (\phi_{i1}^v, ..., \phi_{id}^v)$ and $\phi_i^{pv} = (\phi_{i1}^{pv}, ..., \phi_{id}^{pv})$, where $\phi_{is}^v$ are vector fields and $d$ is a geometric dimension, by solving the following cell problem
\begin{equation}
\label{eq:constraint19}
\begin{split}
\lambda(\cdot,\cdot)^{-1} \phi_{is}^{v} + \nabla \phi_{is}^{pv} = 0,\\
div(\phi_{is}^{v}) = 0,
\end{split}
\end{equation}
where $\int_{R_\omega} \phi_{isk}^{v} \psi_j = \delta_{ij} \delta_{sk} \int_{R_\omega} \psi_j$ and $\int_{R_\omega} \phi_{is}^{pv} \psi_j = 0$.

We can obtain $\phi_i^{vp}$ and $\phi_i^{p}$ by solving the following cell problem
\begin{equation}
\label{eq:constraint20}
\begin{split}
\lambda(\cdot,\cdot)^{-1} \phi_{i}^{vp} + \nabla \phi_{i}^{p} = 0,\\
div(\phi_i^{vp}) = 0,
\end{split}
\end{equation}
where $\int_{R_\omega} \phi_{ik}^{vp} \psi_j = 0$ and $\int_{R_\omega} \phi_{i}^{p} \psi_j = \delta_{ij} \int_{R_\omega} \psi_j$. One can obtain $\phi_{i}^{vp, m}$ and $\phi_{i}^{p, m}$ by solving the same cell problem with $\int_{R_\omega} \phi_{ik}^{vp, m} \psi_j = 0$ and a gradient constraint $\int_{R_\omega} \phi_i^{p, m} \psi_j = \delta_{ij} \int_{R_\omega} (x_m - \tilde{x}_m) \psi_j$.

The multicontinuum basis functions $\phi_i^{vc}$ and $\phi_i^{pc}$ solve the local problem in each cell with no flow boundary conditions and the concentration $\phi_i^c$,
\begin{equation}
\label{eq:constraint22}
\begin{split}
\lambda(\cdot,\cdot)^{-1} \phi_i^{vc} + \nabla\phi^{pc}_i -\phi^{c}_i e_1 =0,\\
div(\phi_i^{vc})=0,\\
\int_{R_\omega^{+}} v \phi^{c}_i\nabla \xi - \sum_{j,l}{\beta_{ij}^l \over \int_{R_\omega^{l}\psi_j^l}}\int_{R_\omega^{l}} \xi \psi_j^l =0,\\
\int_{R_\omega^{l}} \phi^{c}_i\psi_j^l = \delta_{ij}\int_{R_\omega^{l}} \psi_j^l.
\end{split}
\end{equation}
As we mentioned earlier in Section \ref{sec:galerkin_approach}, in our numerical simulations, we take $\phi_i^c$ to be piecewise constant functions, i.e., $\phi_i^c = \chi_i \psi_i$.

Then, we can define the following approximations of the trial and test functions in a vector form
\begin{equation}
M = \begin{bmatrix}
\phi_i^v V_i + \phi_i^{vp} P_i + \phi_i^{vp,m} \nabla_m P_i + \phi_i^{vc} C_i\\
\phi_i^{pv} V_i + \phi_i^p P_i + \phi_i^{p,m} \nabla_m P_i\\
\phi_i^c C_i
\end{bmatrix}, \quad
S_v = \begin{bmatrix}
\phi_j^v W_j\\
\phi_j^{pv} W_j\\
\end{bmatrix}, \quad
S_p = \begin{bmatrix}
\phi_j^{vp} Q_j + \phi_j^{vp,n} \nabla_n Q_j\\
\phi_j^p Q_j + \phi_j^{p,n} \nabla_n Q_j\\
\end{bmatrix}.
\end{equation}

Next, in $R_\omega$, we obtain the following approximation of the variational formulation of the multicontinuum pressure equations
\begin{equation}
\begin{split}
\int_{R_\omega} S_p^T A M =
\int_{R_\omega}\{ (\phi_j^{vp} Q_j + \phi_j^{vp,n} \nabla_n Q_j) \lambda^{-1} (\phi_i^v V_i + \phi_i^{vp} P_i + \phi_i^{vp,m} \nabla_m P_i + \phi_i^{vc} C_i) +\\
(\phi_j^{vp} Q_j + \phi_j^{vp,n} \nabla_n Q_j) \nabla (\phi_i^{pv} V_i + \phi_i^p P_i + \phi_i^{p,m} \nabla_m P_i) -\\
(\phi_j^{vp} Q_j + \phi_j^{vp,n} \nabla_n Q_j) e_1 \phi_i^c C_i +\\
(\phi_j^p Q_j + \phi_j^{p,n} \nabla_n Q_j) div(\phi_i^v V_i + \phi_i^{vp} P_i + \phi_i^{vp,m} \nabla_m P_i + \phi_i^{vc} C_i)\} =\\
\int_{R_\omega} f (\phi_j^p Q_j + \phi_j^{p,n} \nabla_n Q_j).
\end{split}
\end{equation}

Therefore, we have the following multicontinuum pressure equations in the global form
\begin{equation}
\begin{split}
\alpha_{ij}^p V_j 
- \nabla_n (\alpha_{ijn}^p V_j) 
+ \beta_{ij}^p P_j 
- \nabla_n(\overline{\beta}_{ijn}^p P_j)
+ \beta_{ijm}^p \nabla_m P_j
- \nabla_n(\beta_{ijmn}^p \nabla_m P_j)
+ \gamma_{ij}^p C_j
- \nabla_n(\overline{\gamma}_{ijn}^p C_j) = 
f_j^p.
\end{split}
\end{equation}

One can obtain the variational formulation of the multicontinuum velocities in a similar way.
\begin{equation}
\begin{split}
\int_{R_\omega} S_v^T A M =
\int_{R_\omega}\{ \phi_j^{v} W_j \lambda^{-1} (\phi_i^v V_i + \phi_i^{vp} P_i + \phi_i^{vp,m} \nabla_m P_i + \phi_i^{vc} C_i) +\\
\phi_j^{v} W_j \nabla (\phi_i^{pv} V_i + \phi_i^p P_i + \phi_i^{p,m} \nabla_m P_i) -\\
\phi_j^{v} W_j e_1 \phi_i^c C_i +\\
\phi_j^{pv} W_j div(\phi_i^v V_i + \phi_i^{vp} P_i + \phi_i^{vp,m} \nabla_m P_i + \phi_i^{vc} C_i)\} =\\
\int_{R_\omega} f \phi_j^{pv} W_j.
\end{split}
\end{equation}

As a result, we get the following multicontinuum velocity equations in the global form
\begin{equation}
\begin{split}
\alpha_{ij}^v V_j + \beta_{ij}^v P_j + \beta_{ijm}^v \nabla_m P_j + \gamma_{ij}^v C_j = f_j^v.
\end{split}
\end{equation}

The obtained mixed multicontinuum flow model is based on the general multicontinuum expansions and represents the general case of the multicontinuum equations. However, one can consider more simplified cell problems and multicontinuum expansions.

\section{Numerical results}\label{sec:numerical_results}


This section presents the numerical results of our proposed dynamic multicontinuum approach. We consider cases in which continua are characterized by highly contrasted properties as well as gravitational effects. We consider model problems defined in two- and three-continuum media using both mixed and Galerkin multicontinuum modeling approaches.


Let us introduce a computational domain $\Omega = [0, L_1] \times [0, L_2]$, where $L_1 = 9$ and $L_2 = 3$. We will refer to this domain as the target domain, which is of primary interest to us. Multicontinuum modeling on a coarse grid will be performed in the target domain. However, in most of the problems, fine-scale computations will be performed on extended domains: (a) the two-sided extended domain $\Omega_{\text{ext}}^{\text{(LR)}} = [-\Delta^{\text{ext}}, L_1 + \Delta^{\text{ext}}] \times [0, L_2]$ or (b) the right-extended domain $\Omega_{\text{ext}}^{\text{(R)}} = [0, L_1 + 2 \Delta^{\text{ext}}] \times [0, L_2]$, where $\Delta^{\text{ext}} = 1.8$. This is intended to ease the computations due to the dynamic nature of the problems, where coarse blocks and edges may have different numbers of continua at some time instances. The target domain $\Omega$ assures the presence of all the continua. The introduction of the extended domains is not a mandatory requirement of our approach, but rather a way to simplify computations, e.g., one of the model problems is solved without it.

\begin{figure}[hbt!]
\centering
\includegraphics[width=0.49\textwidth]{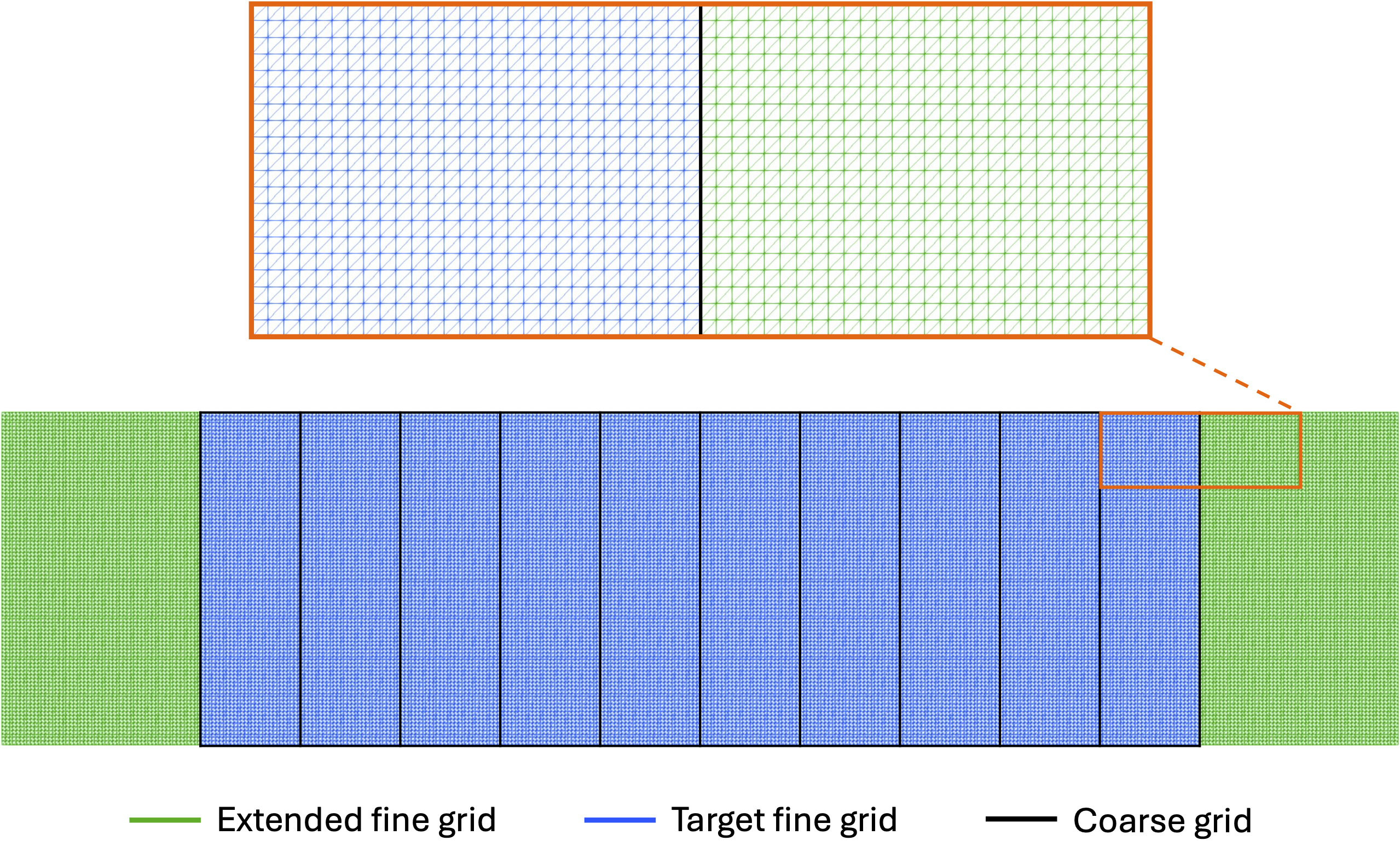}
\includegraphics[width=0.49\textwidth]{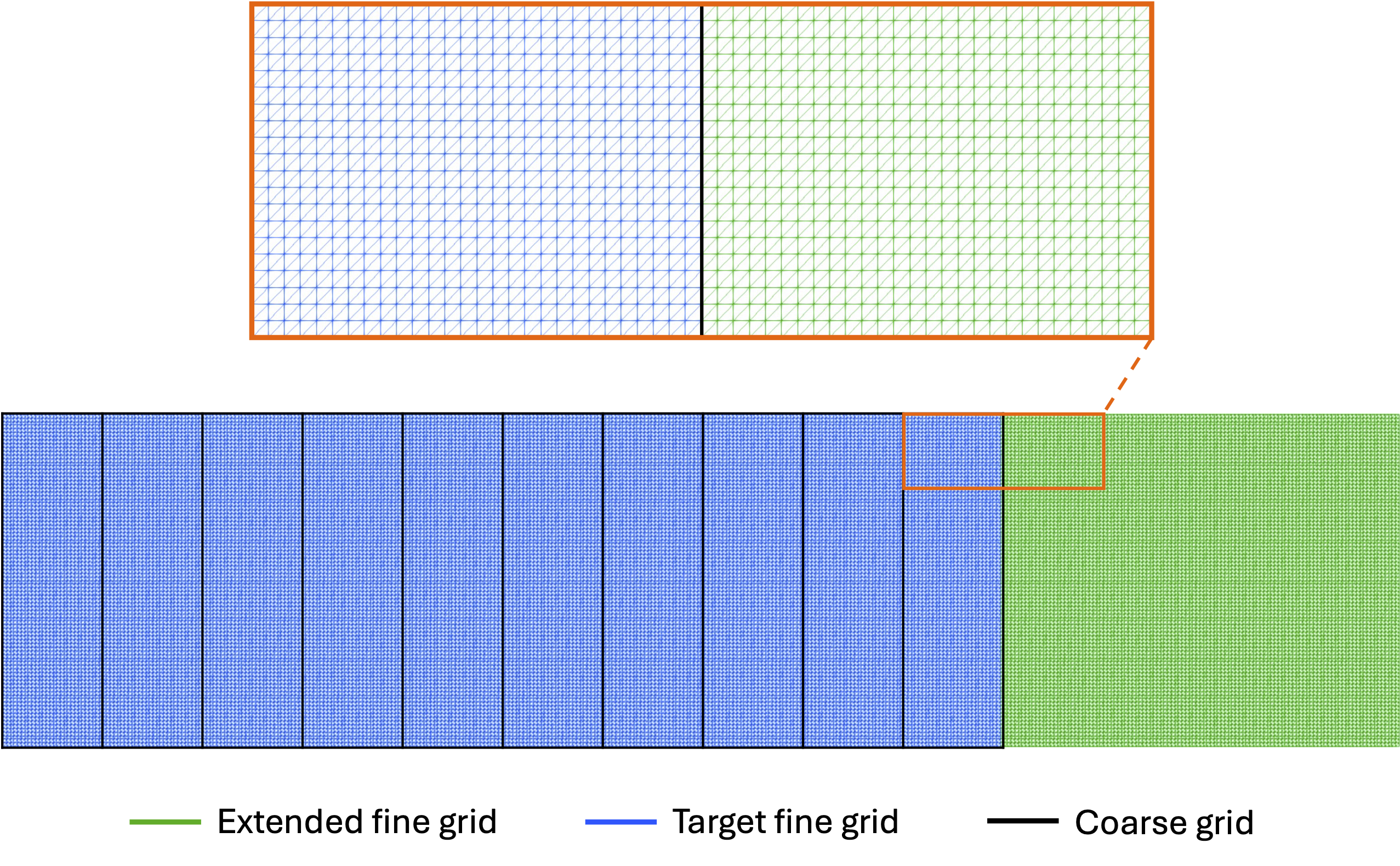}
\caption{Computational grids in the cases of the two-sided extended domain and the right-extended domain (from left to right).}
\label{fig:grids}
\end{figure}

For both cases of the extended domains, the same fine-scale discretizations are used: $\Omega_{\text{ext}}^{\text{(LR)}}$ and $\Omega_{\text{ext}}^{\text{(R)}}$ are partitioned into $392 \times 90$ rectangles, each of which is divided into two triangles. 
The fine grid of the target domain $\Omega$ is the restriction of the extended-domain fine grids to the target domain; it partitions 
$\Omega$ into $280 \times 90$ rectangles, each divided into two triangles. Meanwhile, the coarse grid partitions $\Omega$ into $10 \times 1$ rectangular blocks. Note that, for simplicity, each whole coarse block will be taken as an RVE for itself. Therefore, the terms RVEs and coarse blocks can be used interchangeable in the numerical results. Figure \ref{fig:grids} provides an illustration of the computational grids.

We define continua based on the fine-scale concentration field. In the case of dual-continuum media, we define $I_1 = \{0.5 \leq c \leq 1 \}$ and $I_2 = \{0 \leq c < 0.5 \}$. We introduce $I_1 = \{0.8 \leq c \leq 1 \}$, $I_2 = \{0.4 \leq c < 0.8 \}$, and $I_3 = \{0 \leq c < 0.4 \}$ for the triple-continuum media. For each $I_k$, we introduce a characteristic function $\psi_k(c)$, identifying the continuum at each time.


Note that our proposed multicontinuum approach enables us to formulate macroscopic variables and derive new macroscopic models for nonlinear problems, where coefficients depend on fine-grid functions (in our case, the concentration field). Since our primary objective in this work is to show that the proposed macroscopic models can provide an accurate representation of the coarse-grid solution, we use the fine-grid concentration field to compute macroscopic coefficients. In our next works, we will study computations of these macroscopic coefficients based on the coarse-grid concentration field.


In all the model problems, we first solve the flow problem and then the transport problem using the computed velocity. The temporal approximation of the macroscopic concentration model is based on the Forward Euler method. The time approximation of the fine-scale model varies depending on the specific problem. Numerical implementation is based on the FEniCS computational package \cite{logg2012automated}. For the Lagrangian-Eulerian particle-mesh method, we apply the LeoPart library \cite{maljaars2021leopart}. Visualization of the obtained results is performed using the ParaView software \cite{ayachit2015paraview}.


The following subsections consider various cases of dynamic multicontinuum problems. First, we study problems with gravity-determined continua, and then contrast-determined continua.

\subsection{Gravity-determined continua}\label{sec:gravity_determined}


In this subsection, we consider the cases where continua are determined by gravitational effects. The model problems are well-known gravity-driven fingering problems. We consider problems with two and three continua, as well as with constant and heterogeneous coefficients.


Let us consider fine-scale problem formulation. All the problems are defined on the two-sided extended domain $\Omega_{\text{ext}}^{\text{(LR)}}$ (see Figure \ref{fig:grids}) and have the following general form
\begin{equation}\label{eq:fine_problem_gravity_driven_fingering}
\begin{split}
v= -\lambda(x) (\nabla p-c e_1) \quad &\text{in } \Omega_{\text{ext}}^{\text{(LR)}},\\
div(v)=0 \quad &\text{in } \Omega_{\text{ext}}^{\text{(LR)}},\\
c_t + \nabla \cdot(v c) = 0 \quad &\text{in } \Omega_{\text{ext}}^{\text{(LR)}} \times (0, T].
\end{split}
\end{equation}

We complement the system \eqref{eq:fine_problem_gravity_driven_fingering} with the impermeable boundary conditions $v \cdot n = 0$ on $\partial \Omega_{\text{ext}}^{\text{(LR)}}$, where $n$ is a unit outward normal vector. The initial condition for concentration varies depending on the model problem.


For the numerical solution of the flow problem on a fine grid, we use a mixed finite element method with the lowest-order Raviart-Thomas (RT) basis functions for the velocity and piecewise constant basis functions for the pressure. Whereas the transport problem requires special treatment, since the accuracy of the multicontinuum modeling of the gravity-driven fingering process depends on numerical diffusion on a fine scale. For this reason, we use the Lagrangian-Eulerian particle-mesh method based on the LeoPart library \cite{maljaars2021leopart}. This approach allows us to solve pure convection transport problems without introducing numerical diffusion.

In all the model problems, we randomly seed 2257920 particles. We perform particles advection using the three-stage Runge-Kutta scheme. Note that, in each time step, we first solve the flow problem and then the transport problem using the computed velocity. For particles-mesh projection, we apply a PDE-constrained particle-mesh projection method, developed in \cite{maljaars2019conservative}. To prevent over- and undershoots, we set the penalizing gradients parameter $\zeta$ to a value of 32, which is approximately equal to the average particles count per cell.


In multicontinuum simulations, we take $\phi_i^c = \chi_i \psi_i$ for the transport problem, where $\chi_i = 1 / \int_{R_\omega} \psi_i$, and use a mixed multicontinuum modeling approach for the flow problem. Since $\lambda$ does not possess high contrast in our gravity-driven fingering problems (high contrast cases are considered later), we suppose that we have a single-continuum pressure and multicontinuum velocity fields. Therefore, we use piecewise constant functions as pressure basis functions. While the velocity basis functions are computed by solving simplified cell problems.

The first cell problem is defined in a local domain $\omega_l$, which is constructed as a union of coarse blocks sharing the coarse-grid edge $E_l$. This cell problem accounts for average velocities in different continua and has the following form
\begin{equation}
\begin{split}
\lambda^{-1}(x) \phi_{i}^{v} + \nabla \phi_{i}^{pv} = 0 \quad \text{in } \omega_l,\\
div (\phi_{i}^{v}) = \sigma_{i}^{v, l} \quad \text{in } \omega_l,
\end{split}
\end{equation}
with a boundary condition $\phi_i^{v} \cdot n = 0$ on $\partial \omega$ and an additional constraint $\phi_{i}^{v} \cdot n = \psi_{i}(c)$ on $E_l$. Note that $\sigma_{i}^{v, l}$ is chosen to satisfy the compatibility condition $\int_{K} \sigma_{i}^{v, l} = \int_{E_l} \psi_i(c)$ for all $K \subset \omega_l$.

The second cell problem accounts for the gravitational effects and is formulated in each coarse block separately. One can represent this cell problem in the following way
\begin{equation}
\begin{split}
\lambda^{-1}(x) \phi_{i}^{vc} + \nabla \phi_{i}^{pc} = \psi_{i}(c) e_1 \quad \text{in } K_l,\\
div (\phi_{i}^{vc}) = 0 \quad \text{in } K_l.
\end{split}
\end{equation}
Here, we set a boundary condition $\phi_i^{vc} \cdot n = 0$ on $\partial K_l$.


Next, we present various cases of gravity-driven fingering problems, including two and three continua, as well as considering heterogeneous coefficient $\lambda$.

\subsubsection{Gravity-driven fingering: two continua}\label{sec:gravity_fingering_dual_continuum}


This case represents the gravity-driven fingering problem in a dual-continuum medium ($N = 2$). We start with some initial configuration of fingers, which then dynamically evolves during the simulation. The coefficient $\lambda$ is set to 1. We simulate for 120 time steps with a step size of $\tau = 3.33 \cdot 10^{-2}$, the total time is $T = 3.996$. The Courant-Friedrichs-Lewy (CFL) number of the fine-scale simulation varies from 0.442 to 0.487 over time. The CFL number is defined as $\text{CFL}(t) = \dfrac{\tau}{h} \max_{x \in \Omega_{\text{ext}}^{\text{(LR)}}} |v(x, t)|$, where $h$ is the fine-grid cell diameter. The cases in Section \ref{sec:contrast_determined} use the same definition with their respective domains, i.e., $\Omega_{\text{ext}}^{\text{(R)}}$ or $\Omega$.

\begin{figure}[hbt!]
\centering
\includegraphics[width=0.49\textwidth]{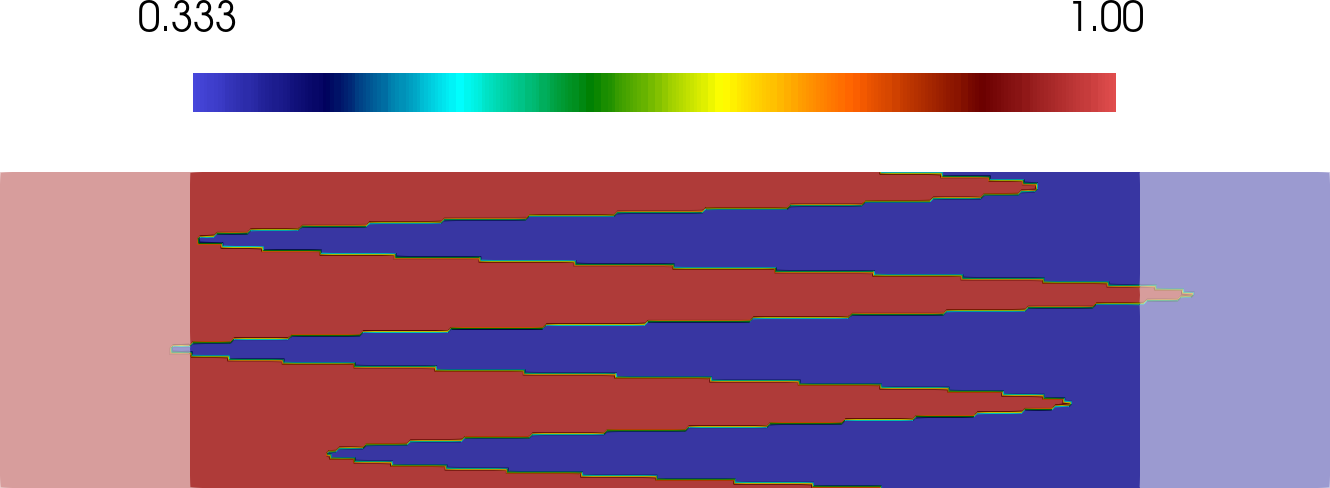}
\includegraphics[width=0.49\textwidth]{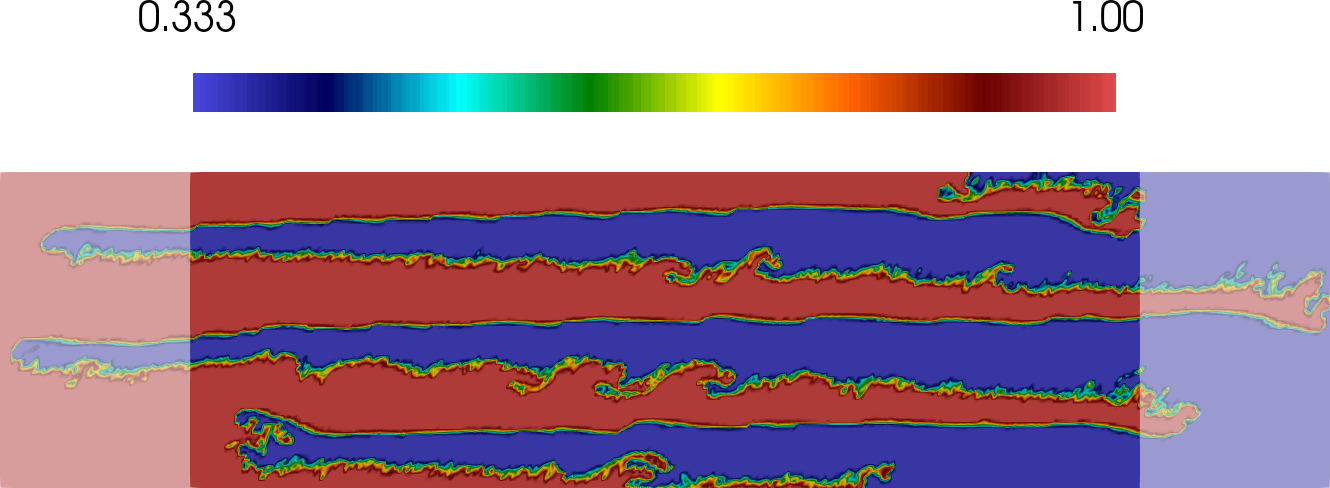}
\caption{Initial and final fine-scale concentration distributions. Gravity-driven fingering with two continua.}
\label{fig:c_fine_problem_3_N_2}
\end{figure}

Figure \ref{fig:c_fine_problem_3_N_2} depicts the initial and final time concentration distributions (from left to right) on the fine grid. The extension regions on the sides are indicated by semi-transparent shading. We can see that the configuration of the fingers has changed significantly during the simulation. The first continuum fingers, characterized by a higher concentration and representing a denser fluid, are directed to the right in accordance with the direction of gravity. Meanwhile, the second continuum fingers, corresponding to a less dense fluid, are directed in the opposite direction. At the same time, we see the development of smaller fingers.

\begin{figure}[hbt!]
\centering
\includegraphics[width=0.49\textwidth]{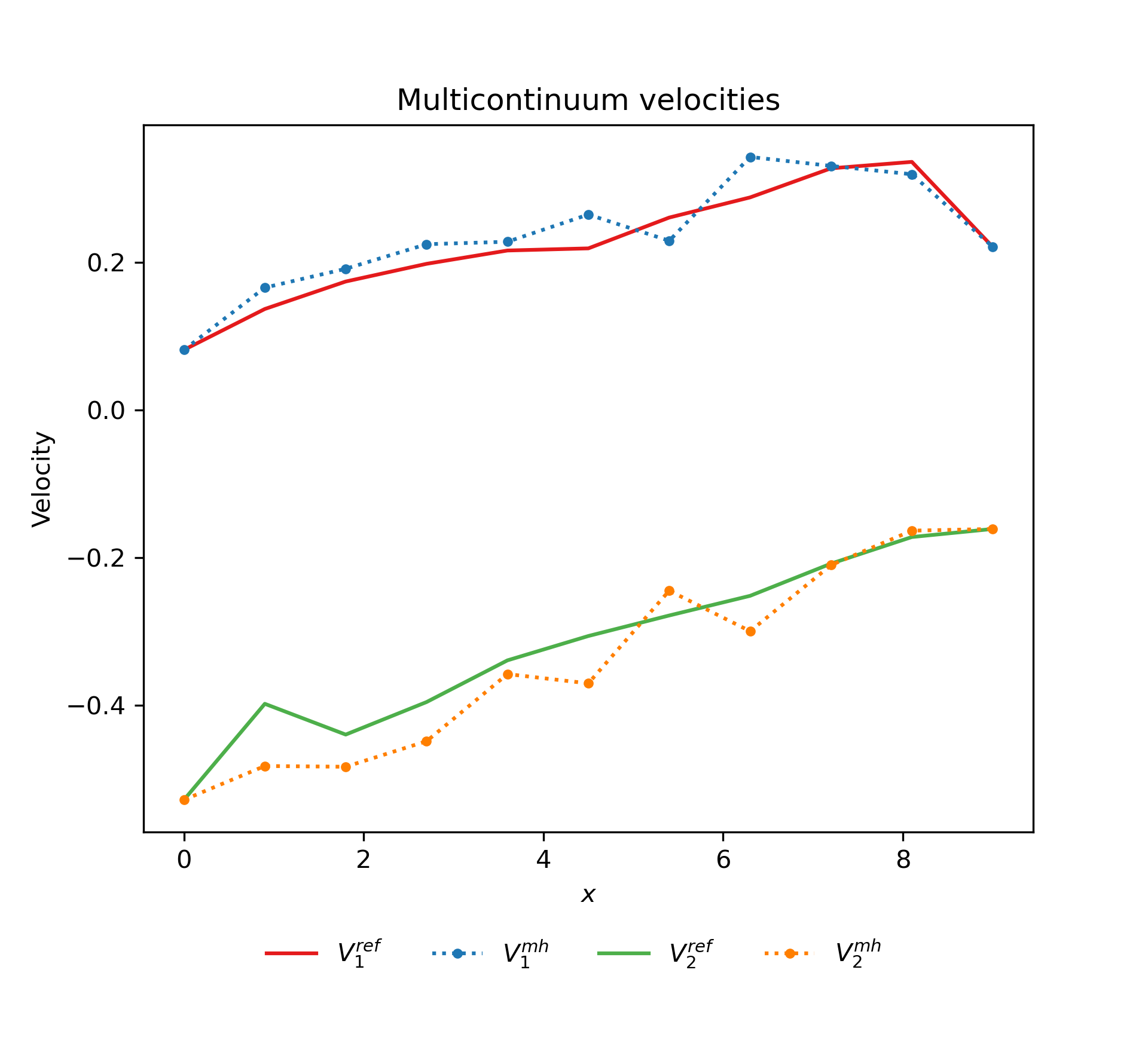}
\includegraphics[width=0.49\textwidth]{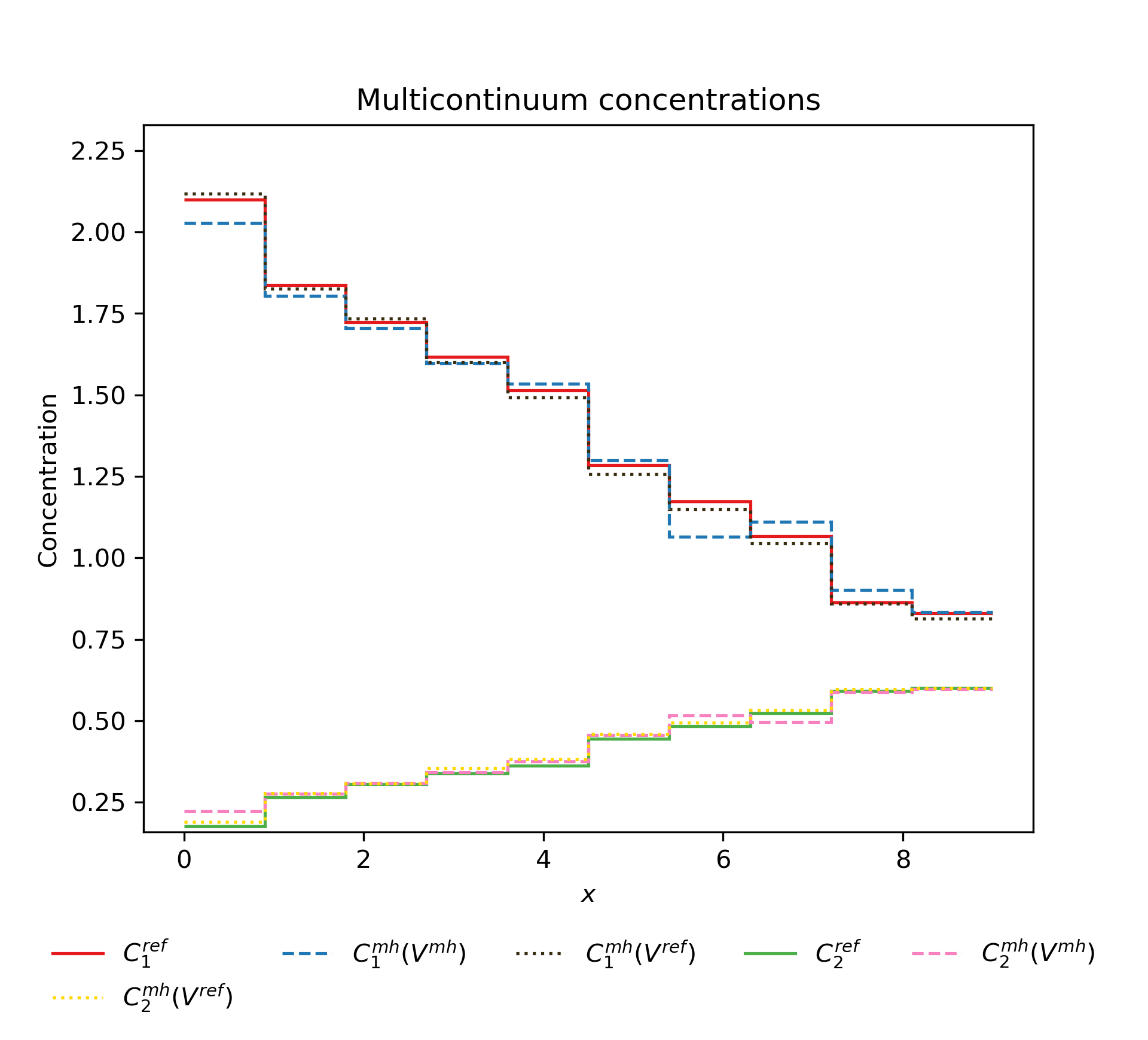}
\caption{Multicontinuum velocity and concentration plots (from left to right) at the final time. Gravity-driven fingering with two continua.}
\label{fig:coarse_results_problem_3_N_2}
\end{figure}

In Figure \ref{fig:coarse_results_problem_3_N_2}, we present plots of multicontinuum velocities and concentrations at the final time (from left to right). These plots depict both the reference average solutions and the solutions obtained using the proposed multicontinuum approach. One can see from the left figure that our homogenized velocities $V_i^{\text{mh}}$ accurately approximate the reference average velocities $V_i^{\text{ref}}$. The right figure depicts plots of the reference multicontinuum concentrations $C_i^{\text{ref}}$ and two types of homogenized concentrations: $C_i^{\text{mh}} (V^{\text{ref}})$ and $C_i^{\text{mh}} (V^{\text{mh}})$. In $C_i^{\text{mh}}(V^{\text{ref}})$, we use the reference average velocities $V^{\text{ref}} = (V_1^{\text{ref}}, V_2^{\text{ref}})$ to compute the homogenized concentrations $C_i^{\text{mh}}$, and in $C_i^{\text{mh}}(V^{\text{mh}})$, we use the homogenized velocities $V^{\text{mh}} = (V_1^{\text{mh}}, V_2^{\text{mh}})$. We can see that all the concentration plots match, indicating the high accuracy of our model.

One can see that the concentration of the first continuum exceeds 1. This is because we define our multicontinuum concentrations as $C_i (x_\omega) = \int_{R_\omega} c \psi_i(c)$ (see \eqref{eq:macroscopic_variable_definitions}) instead of $\bar{C}_i (x_\omega) = \int_{R_\omega} c \psi_i(c) / \int_{R_\omega} \psi_i(c)$. We prefer this type of multicontinuum concentrations, since it is more informative in terms of the continuum dynamics' representation. Moreover, one can obtain $\bar{C}_i$ from the corresponding scaling of $C_i$. 

\begin{figure}[hbt!]
\centering
\begin{subfigure}[b]{\textwidth}
\centering
\includegraphics[width=0.32\textwidth]{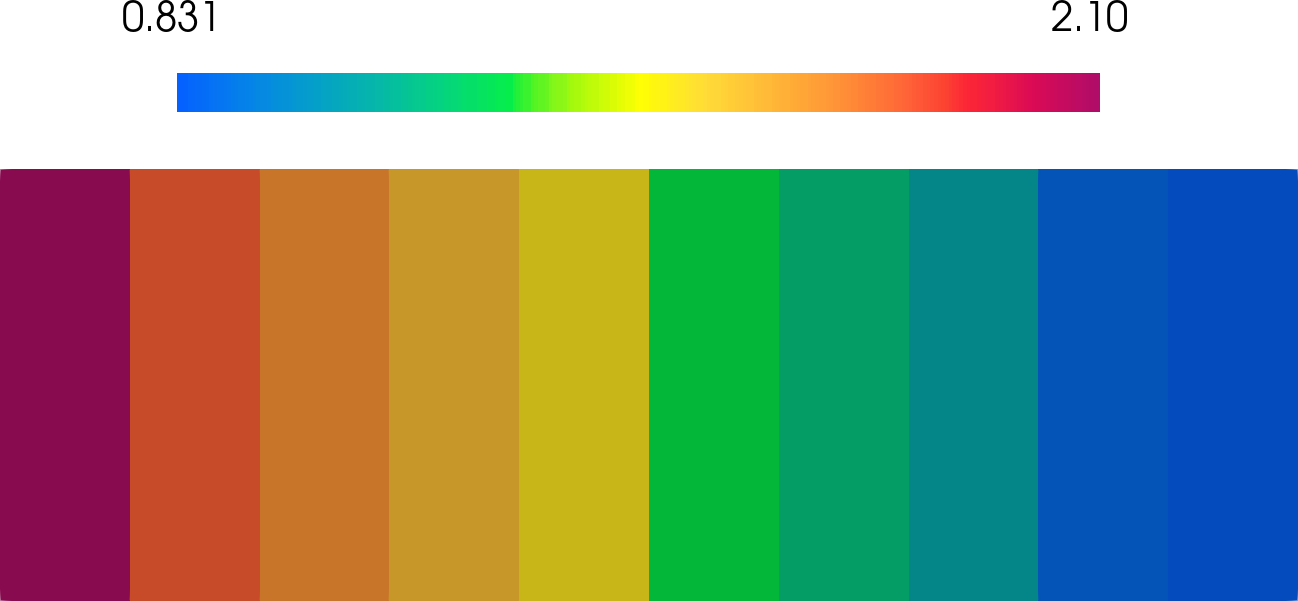}
\includegraphics[width=0.32\textwidth]{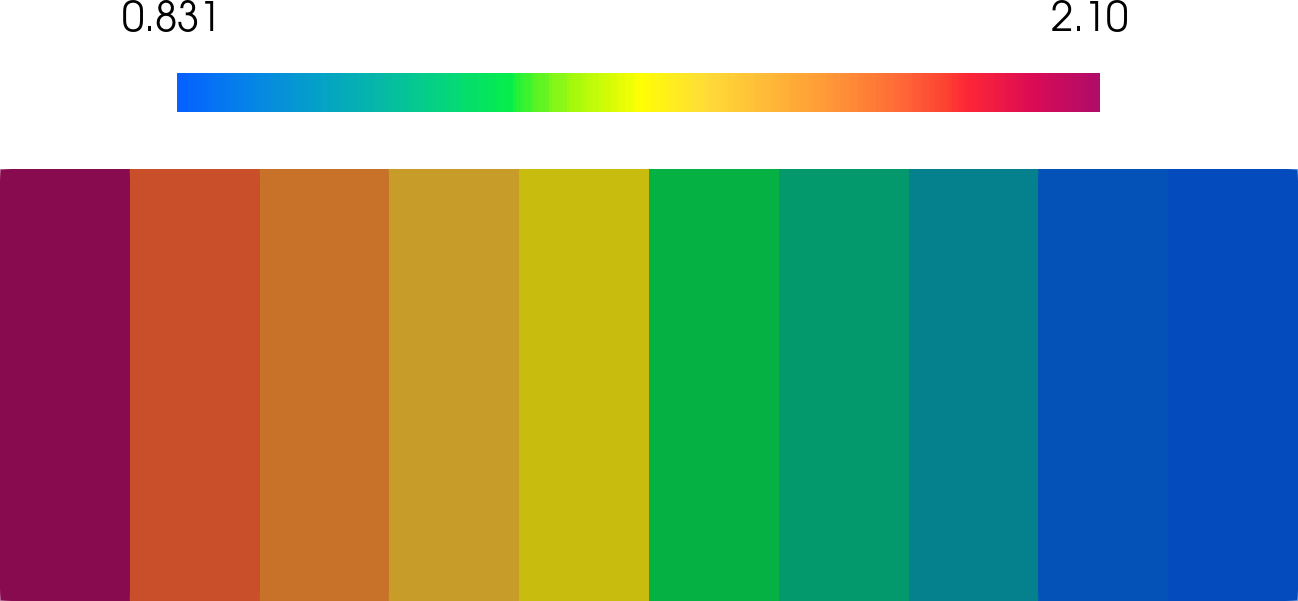}
\includegraphics[width=0.32\textwidth]{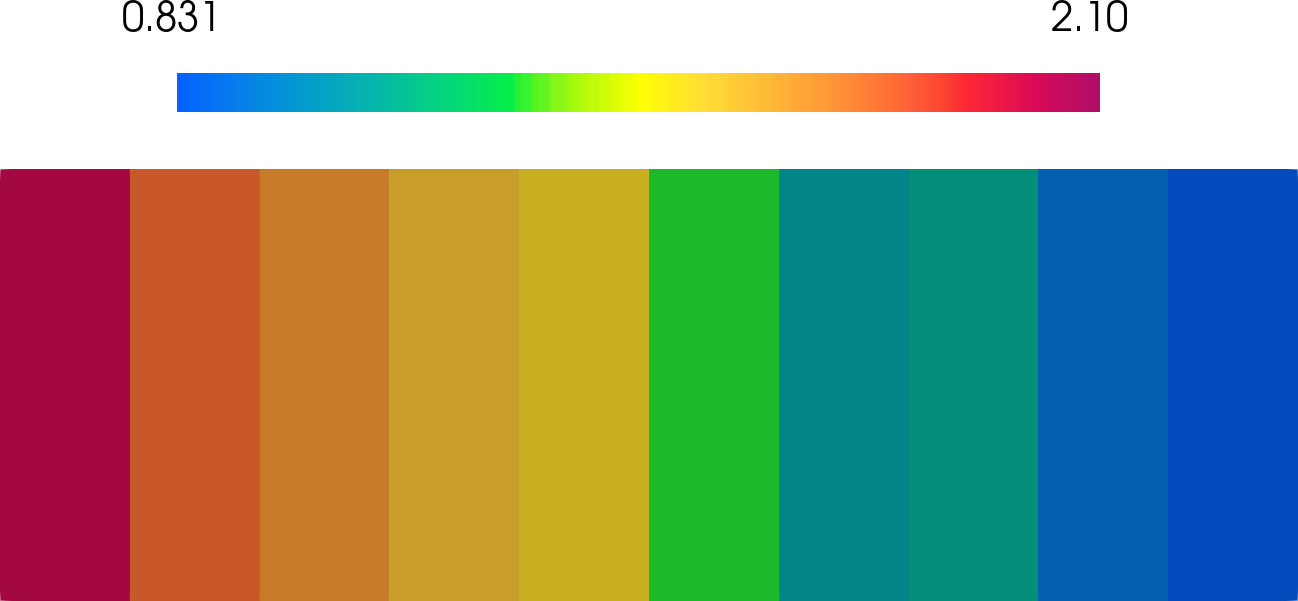}
\caption{First continuum concentrations (left to right): $C_1^{\text{ref}}$, $C_1^{\text{mh}} (V^{\text{ref}})$, and $C_1^{\text{mh}} (V^{\text{mh}})$.}
\end{subfigure}
\begin{subfigure}[b]{\textwidth}
\centering
\includegraphics[width=0.32\textwidth]{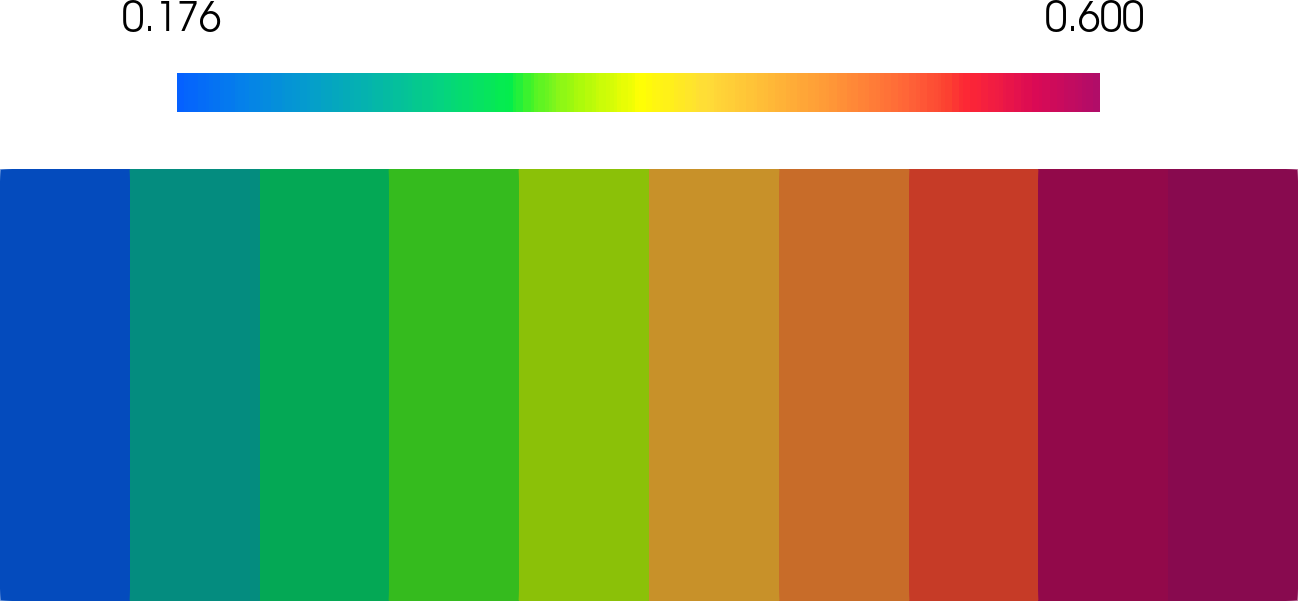}
\includegraphics[width=0.32\textwidth]{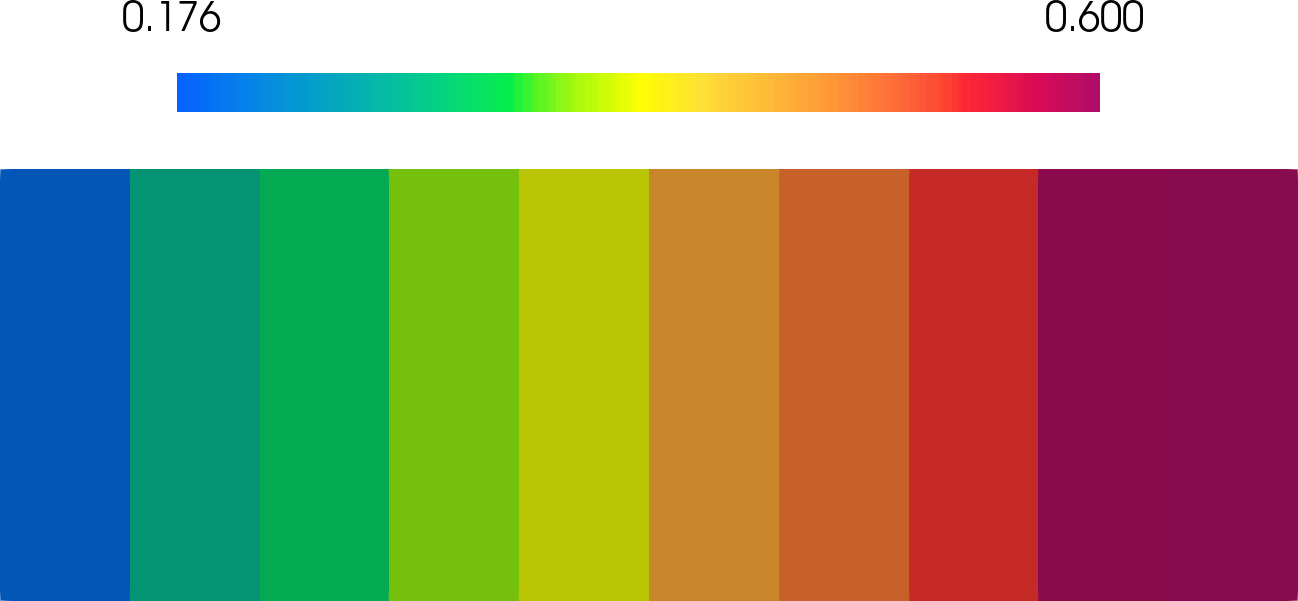}
\includegraphics[width=0.32\textwidth]{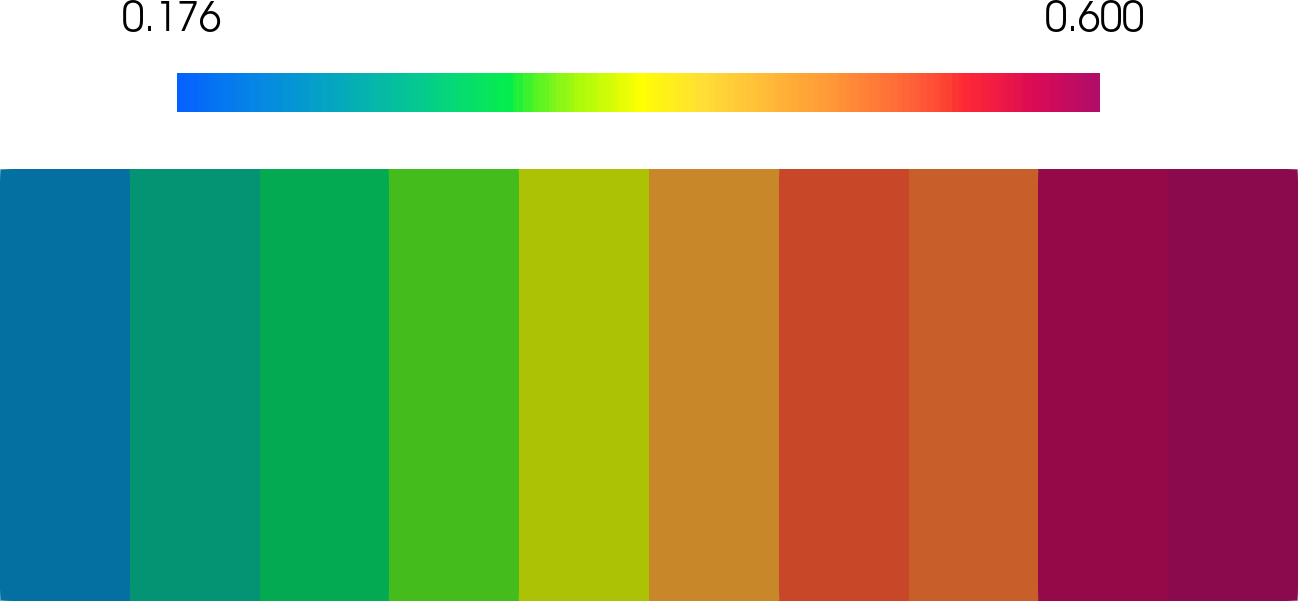}
\caption{Second continuum concentrations (left to right): $C_2^{\text{ref}}$, $C_2^{\text{mh}} (V^{\text{ref}})$, and $C_2^{\text{mh}} (V^{\text{mh}})$.}
\end{subfigure}
\caption{Multicontinuum concentrations at the final time. Gravity-driven fingering with two continua.}
\label{fig:C_color_bar_results_problem_3_N_2}
\end{figure}

Figure \ref{fig:C_color_bar_results_problem_3_N_2} presents two-dimensional representations of the multicontinuum concentrations. From top to bottom, we depict the first and second continuum concentrations. From left to right, we present the reference concentration, the homogenized concentration based on the reference velocities, and the homogenized concentration based on the homogenized velocities. One can see that all the figures are very similar, which indicates the high accuracy of our multicontinuum model.

Let us consider the errors in the obtained multiscale solutions. We compare our solutions with the reference average velocities and concentrations. For the velocities, we use the following relative $L^2$ errors for $i = 1, ..., N$
\begin{equation}\label{eq:errors_U}
e_{V}^{(i)} = \frac{\| V_i^{\text{mh}} - V_i^{\text{ref}} \|_2}{\| V_i^{\text{ref}} \|_2} \times 100\%,
\end{equation}
where $\| V_i \|_2 = \sqrt{\sum_E (V_i^E)^2}$, and $V_i^{\text{ref}}$ is the reference average velocity of the continuum $i$, and $V_i^{\text{mh}}$ is the homogenized velocity of the continuum $i$.

In our numerical results, we have the following relative $L^2$ errors of the homogenized velocities at the final time: $e_{V}^{(1)} = 11.723$\% and $e_{V}^{(2)} = 12.647$\%. In this way, our proposed multicontinuum approach provides good accuracy for the multicontinuum velocities of the gravity-driven fingering process in a dual-continuum medium.

Let us now consider the errors of the multicontinuum concentrations. For a more comprehensive analysis, we calculate several types of relative errors. First, we compare the homogenized concentrations computed using reference velocities $C_i^{\text{mh}} (V^{\text{ref}})$ with the reference concentrations $C_i^{\text{ref}}$. This type represents our irreducible error. Then, we compare the homogenized concentrations computed using homogenized velocities $C_i^{\text{mh}} (V^{\text{mh}})$ with the reference concentrations $C_i^{\text{ref}}$. This type corresponds to the standard error with the reference solution. At last, we compare $C_i^{\text{mh}} (V^{\text{mh}})$ with $C_i^{\text{mh}} (V^{\text{ref}})$ to get an error with a reduced influence of coarse-grid numerical diffusion.

\begin{table}[hbt!]
\caption{Relative errors of multicontinuum concentrations at the final time for gravity-driven fingering with two continua.}
\label{tab:C_errors_problem_3_N_2}
\centering
\begin{tabular}{c | c c }
Error $e_{C}^{(i)}$ & $e_{C}^{(1)}$ & $e_{C}^{(2)}$ \\ \hline
$\|C_i^{\text{mh}} (V^{\text{ref}}) - C_i^{\text{ref}} \|_2 / \| C_i^{\text{ref}} \|_2 \times 100\%$ & 1.265\% & 2.876\% \\ 
$\|C_i^{\text{mh}} (V^{\text{mh}}) - C_i^{\text{ref}} \|_2 / \| C_i^{\text{ref}} \|_2 \times 100\%$ & 3.235\% & 4.826\% \\ 
$\|C_i^{\text{mh}} (V^{\text{mh}}) - C_i^{\text{mh}} (V^{\text{ref}}) \|_2  / \| C_i^{\text{mh}} (V^{\text{ref}}) \|_2 \times 100\%$ & 3.552\% & 4.027\% \\ 
\end{tabular}
\end{table}

Table \ref{tab:C_errors_problem_3_N_2} presents different relative errors of the multicontinuum concentrations at the final time. Note that here we have $\| C_i \|_2 = \sqrt{\sum_K (C_i^K)^2}$. One can see that the obtained errors are minor. As expected, the lowest errors are when we compute the homogenized concentrations using the reference velocities. However, the errors of the homogenized concentrations based on the homogenized velocities are also minor. The errors between these two homogenized concentrations are small. Therefore, our proposed multicontinuum model can accurately describe the multicontinuum concentration dynamics for the dual-continuum gravity-driven fingering problem.

\subsubsection{Gravity-driven fingering: three continua}\label{sec:gravity_fingering_triple_continuum}


In this case, we consider the process of gravity-driven fingering in a triple-continuum medium ($N = 3$). We simulate for 150 time steps with a step size of $\tau = 4 \cdot 10^{-2}$, the total time is $T = 6$. We set the coefficient $\lambda = 1$. The CFL number varies from 0.396 to 0.437 over time for fine-grid computations.

\begin{figure}[hbt!]
\centering
\includegraphics[width=0.49\textwidth]{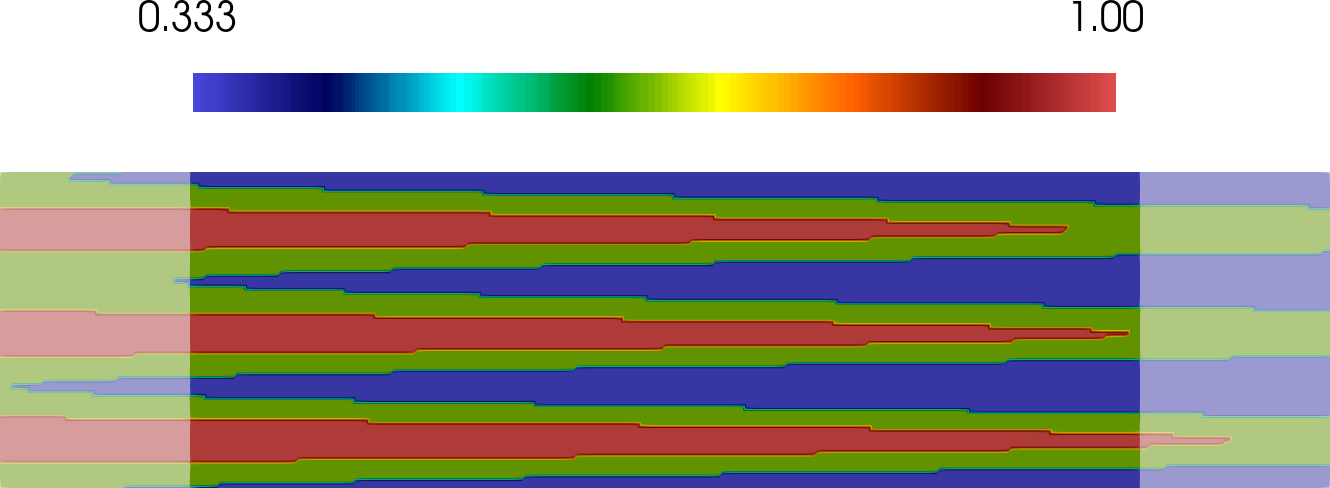}
\includegraphics[width=0.49\textwidth]{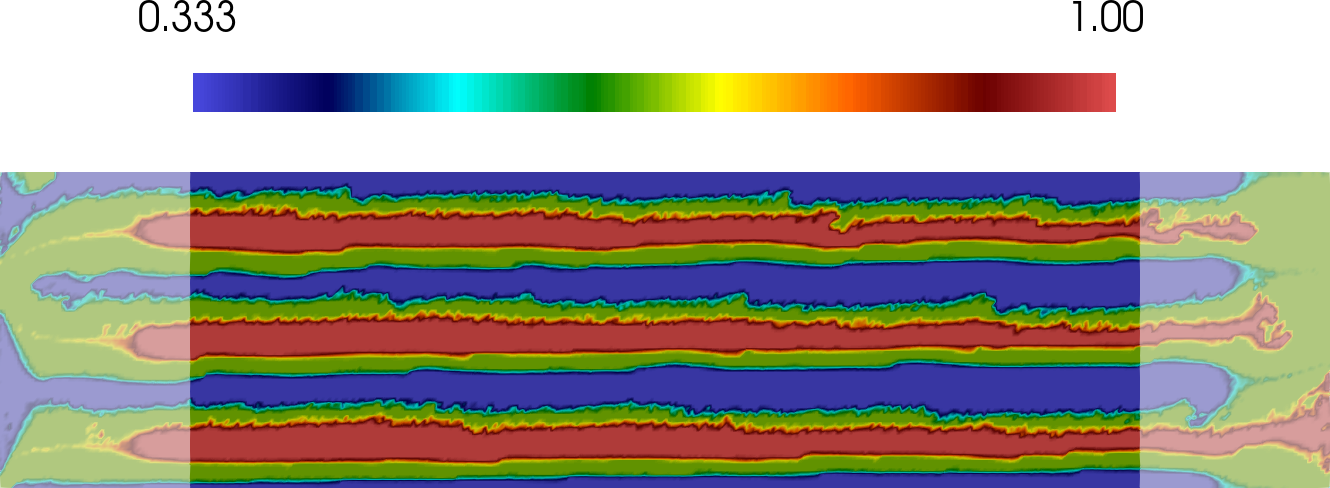}
\caption{Initial and final fine-scale concentration distributions. Gravity-driven fingering with three continua.}
\label{fig:c_fine_problem_3_N_3}
\end{figure}

Figure \ref{fig:c_fine_problem_3_N_3} presents the fine-scale concentration distributions at the initial and final time instances (from left to right). In the initial condition, the first continuum has a concentration of 1, the second continuum has a concentration of 0.666, and the third continuum has a concentration of 0.333. As before, the extended regions are indicated by semi-transparent shading. One can see that the fingers have undergone considerable evolution over time. Their structure has changed, and smaller fingers have developed. The first continuum fingers, corresponding to the densest fluid, are directed to the right in accordance with gravity. The third continuum fingers, representing the least dense fluid, are directed to the left. The second continuum fingers, representing the intermediate-density fluid, experience competing influences from the denser and lighter regions.

\begin{figure}[hbt!]
\centering
\includegraphics[width=0.49\textwidth]{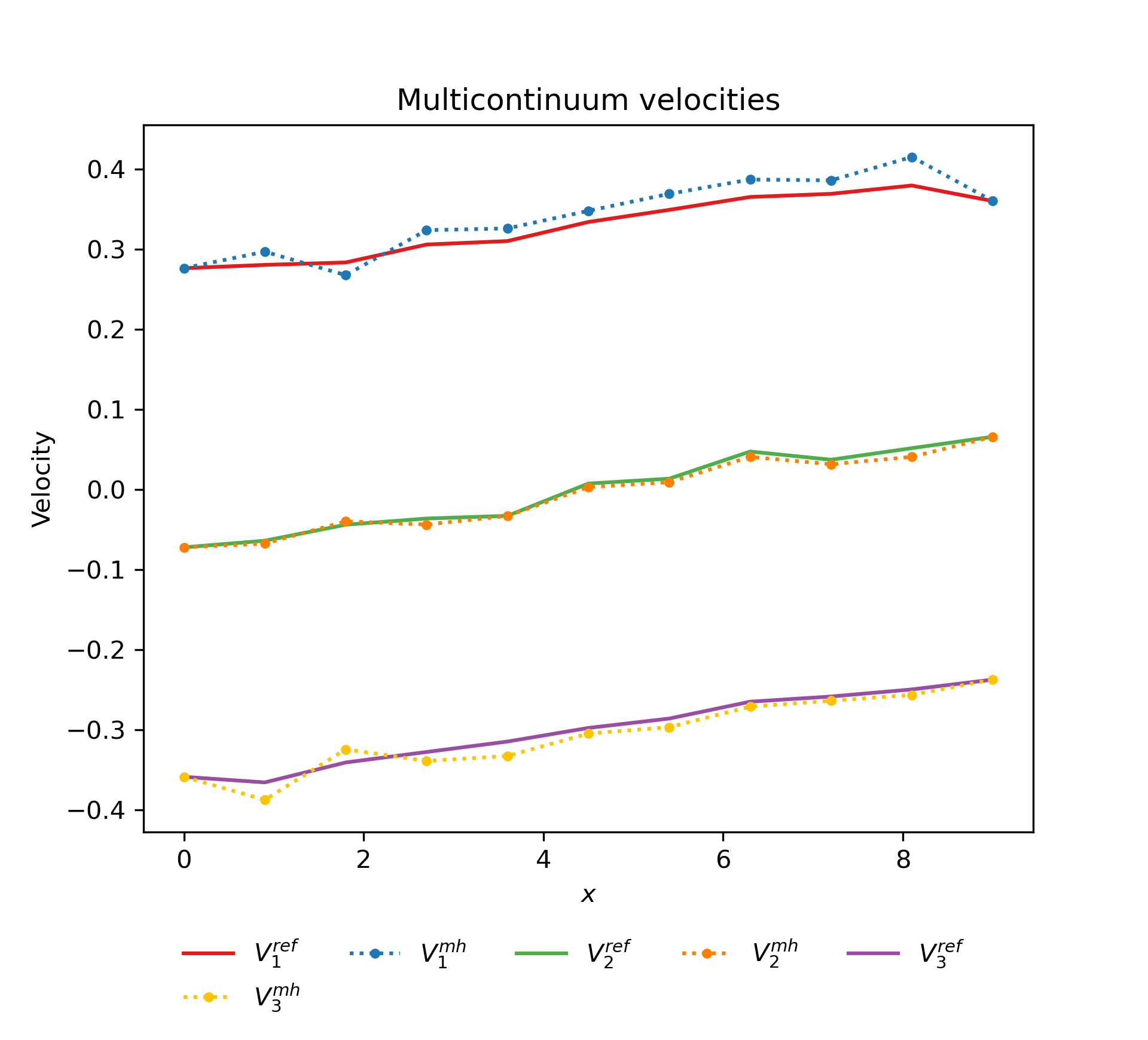}
\includegraphics[width=0.49\textwidth]{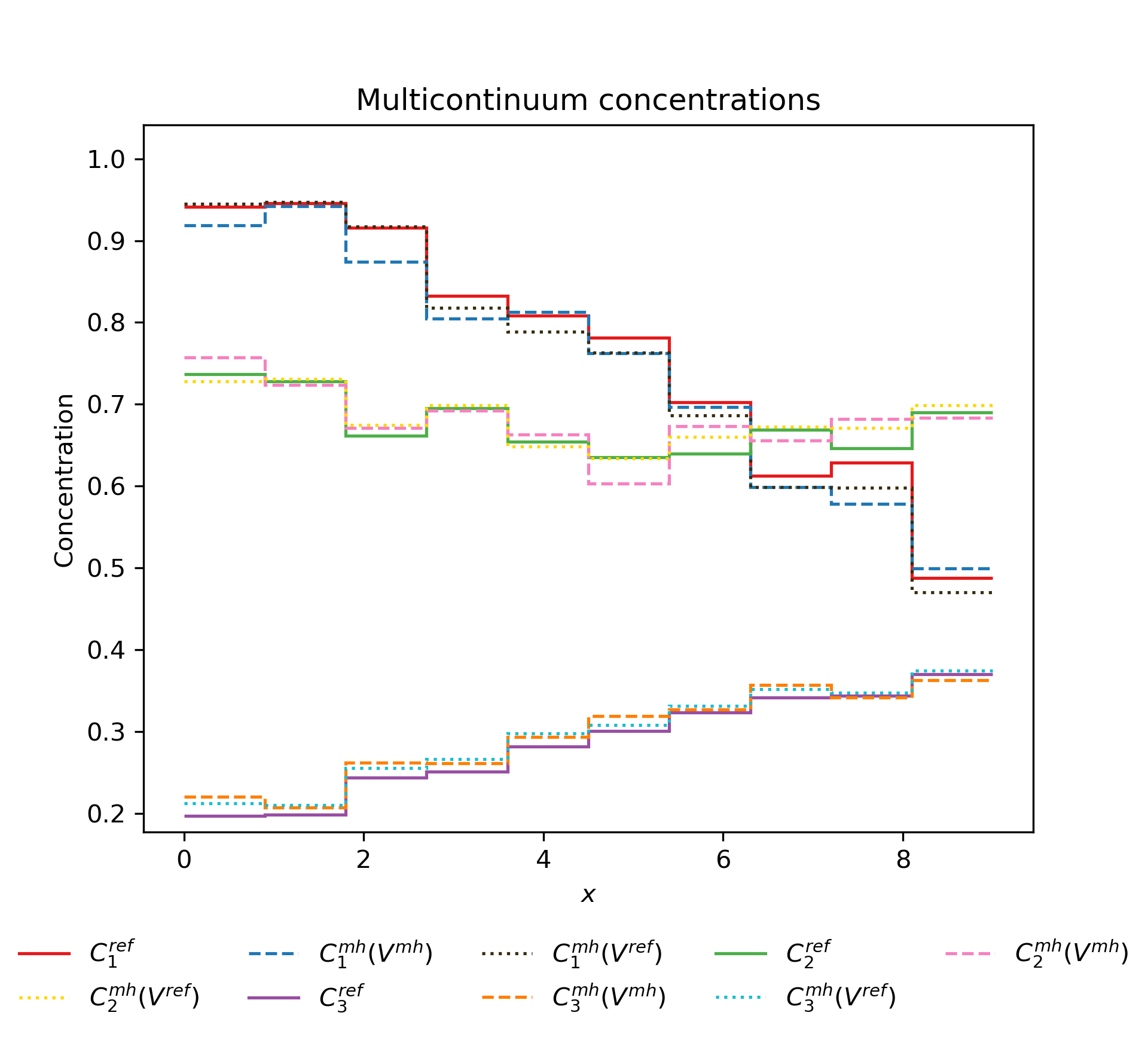}
\caption{Multicontinuum velocity and concentration plots (from left to right) at the final time. Gravity-driven fingering with three continua.}
\label{fig:coarse_results_problem_3_N_3}
\end{figure}

Let us now consider the coarse-scale results. Figure \ref{fig:coarse_results_problem_3_N_3} depicts the multicontinuum velocities and concentrations at the final time (from left to right). We can see that our homogenized velocities approximate the reference multicontinuum velocities with high accuracy. Note that the velocity of the second continuum is close to zero on average. This behavior is caused by the symmetry of the initial concentration distribution (see Figure \ref{fig:c_fine_problem_3_N_3}). In terms of multicontinuum concentrations, we also observe that our multicontinuum solutions demonstrate high accuracy. Both homogenized concentrations (based on $V^{\text{ref}}$ and $V^{\text{mh}}$) accurately follow the trends of the reference concentrations for all three continua.

Next, we present the errors of our multicontinuum solutions. This case requires special treatment for computing velocity errors due to the closeness of the second continuum's velocity to zero. Instead of calculating relative $L^2$ errors \eqref{eq:errors_U}, we consider absolute $L^2$ errors for each continuum. In addition, we compute a general relative $L^2$ error for the whole velocity field $V = (V_1, V_2)$. Therefore, we use the following error norms
\begin{equation}\label{eq:errors_U_N_3}
e_{V, \text{abs}}^{(i)} = \| V_i^{\text{mh}} - V_i^{\text{ref}} \|_2, \quad 
e_{V} = \frac{\| V^{\text{mh}} - V^{\text{ref}} \|_2}{\| V^{\text{ref}} \|_2} \times 100\%.
\end{equation}

As a result, we have the following absolute $L^2$ errors for each continuum at the final time: $e_{V, \text{abs}}^{(1)} = 0.061$, $e_{V, \text{abs}}^{(2)} = 0.018$, $e_{V, \text{abs}}^{(3)} = 0.038$. The general relative $L^2$ error for the whole velocity field is $e_{V} = 4.958$ \%. Therefore, our multicontinuum approach can provide high accuracy for multicontinuum velocities in the triple-continuum case.

\begin{table}[hbt!]
\caption{Relative errors of multicontinuum concentrations at the final time for gravity-driven fingering with three continua.}
\label{tab:C_errors_problem_3_N_3}
\centering
\begin{tabular}{c | c c c }
Error $e_{C}^{(i)}$ & $e_{C}^{(1)}$ & $e_{C}^{(2)}$ & $e_{C}^{(3)}$ \\ \hline
$\|C_i^{\text{mh}} (V^{\text{ref}}) - C_i^{\text{ref}} \|_2 / \| C_i^{\text{ref}} \|_2 \times 100\%$ & 2.090\% & 1.779\% & 3.783\% \\ 
$\|C_i^{\text{mh}} (V^{\text{mh}}) - C_i^{\text{ref}} \|_2 / \| C_i^{\text{ref}} \|_2 \times 100\%$ & 3.208\% & 3.055\% & 4.668\% \\ 
$\|C_i^{\text{mh}} (V^{\text{mh}}) - C_i^{\text{mh}} (V^{\text{ref}}) \|_2  / \| C_i^{\text{mh}} (V^{\text{ref}}) \|_2 \times 100\%$ & 2.803\% & 2.511\% & 2.326\% \\ 
\end{tabular}
\end{table}

Table \ref{tab:C_errors_problem_3_N_3} presents relative $L^2$ errors for the final-time multicontinuum concentrations. As in the previous case, we consider different types of errors to perform a more thorough analysis. One can see that all the errors are minor. We observe that the errors of the homogenized concentrations based on the reference velocities $C_i^{\text{mh}} (V^{\text{mh}})$ are smaller than those for the homogenized concentrations based on the homogenized velocities $C_i^{\text{mh}} (V^{\text{ref}})$, as expected. The errors between these homogenized concentrations, which correspond to the reduction of the numerical diffusion's influence, are also minor. Therefore, our multicontinuum model can accurately approximate the multicontinuum concentrations in the triple-continuum case of gravity-driven fingering.

\subsubsection{Gravity-driven fingering: two continua with heterogeneities}\label{sec:gravity_fingering_heterogeneous_dual_continuum}


Let us now consider the problem of gravity-driven fingering in a heterogeneous dual-continuum medium. Unlike the previous problems, the coefficient $\lambda$ is spatially heterogeneous. The initial condition is set the same as in the dual-continuum case with $\lambda=1$. We set the number of time steps to 120 with a step size of $\tau = 3.33 \cdot 10^{-2}$, the total time is $T = 3.996$. The CFL number varies from 0.278 to 0.297 depending on the time instance.

\begin{figure}[hbt!]
\centering
\includegraphics[width=0.49\textwidth]{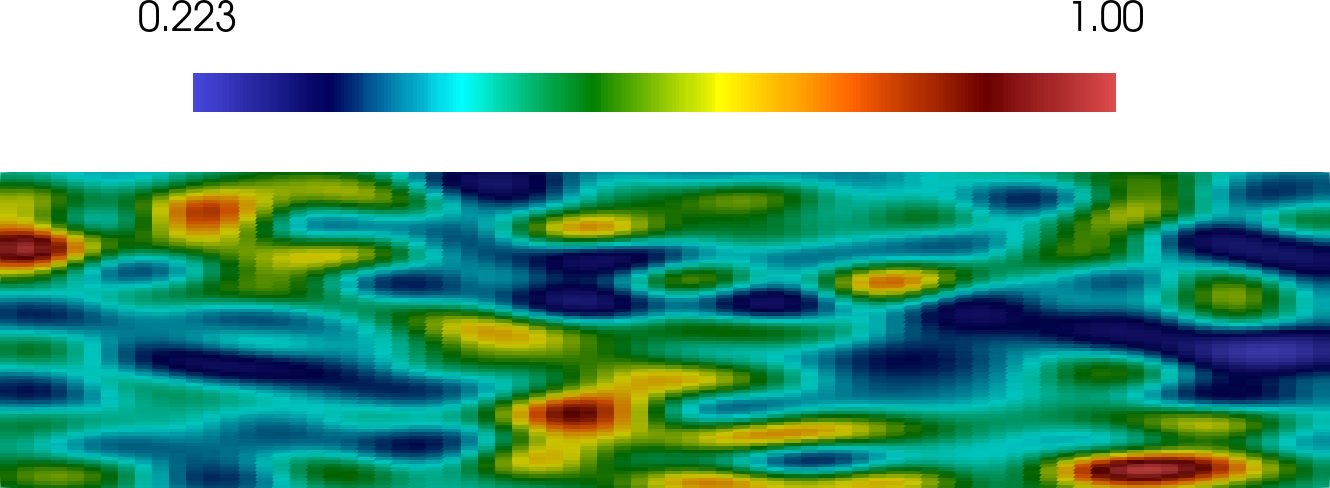}
\includegraphics[width=0.49\textwidth]{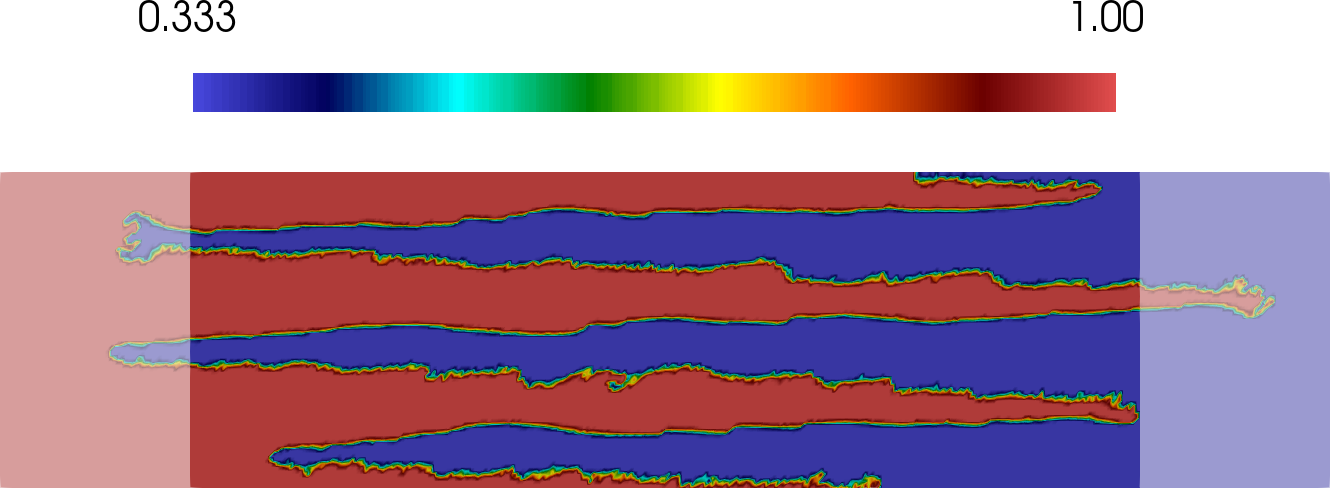}
\caption{Distributions of heterogeneous coefficient $\lambda$ and final fine-scale concentration. Gravity-driven fingering with two continua and heterogeneities.}
\label{fig:c_fine_problem_4_N_2}
\end{figure}

In Figure \ref{fig:c_fine_problem_4_N_2}, we depict the heterogeneous coefficient $\lambda$ and the concentration distribution at the final time. One can clearly see the influence of the heterogeneous coefficient on the evolution of the fingers. Again, the first continuum fingers are directed to the right in accordance with gravity. Whereas the second continuum fingers, corresponding to less dense fluid, are oriented to the opposite side. Smaller secondary fingers are also observed.

\begin{figure}[hbt!]
\centering
\includegraphics[width=0.49\textwidth]{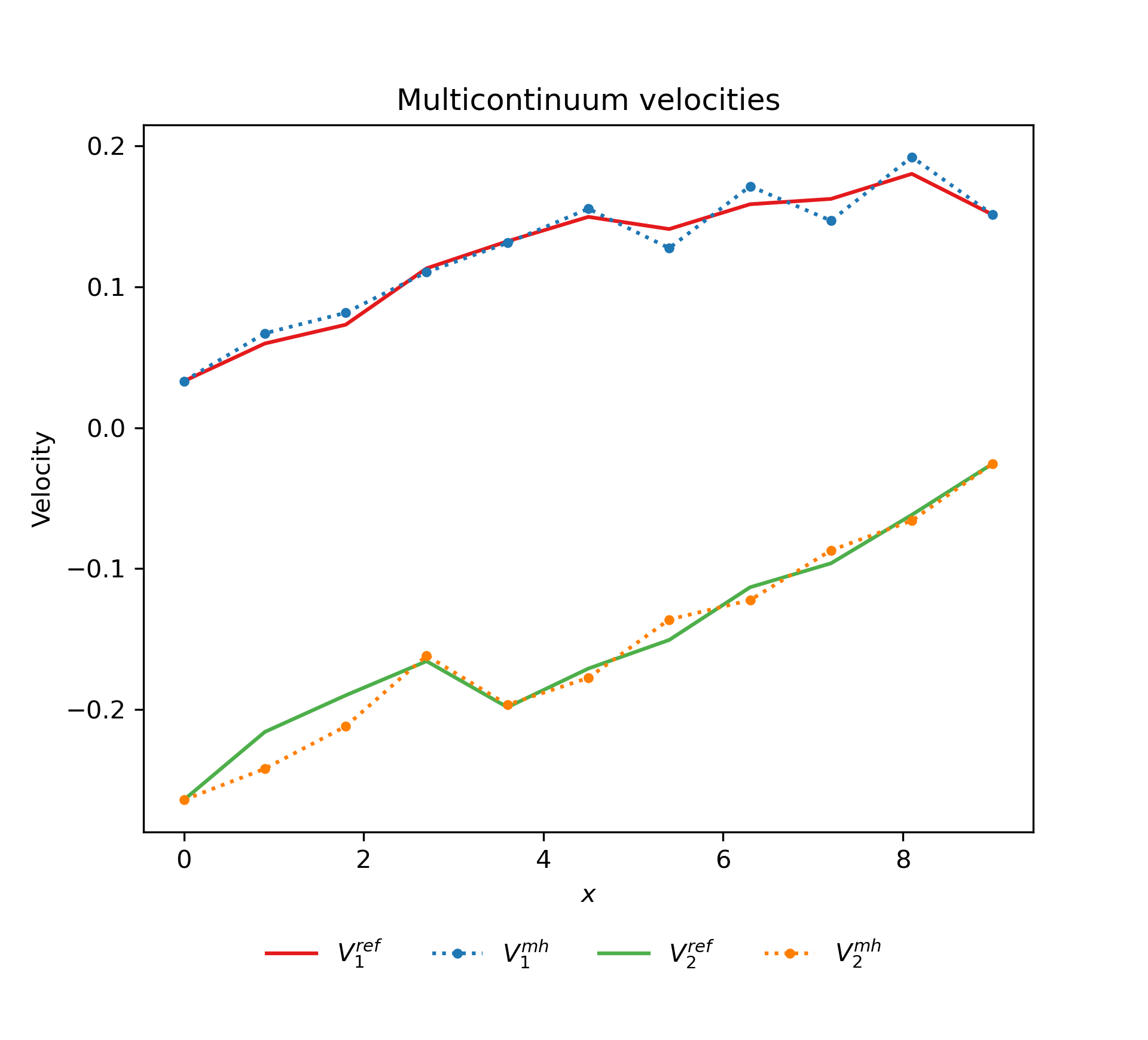}
\includegraphics[width=0.49\textwidth]{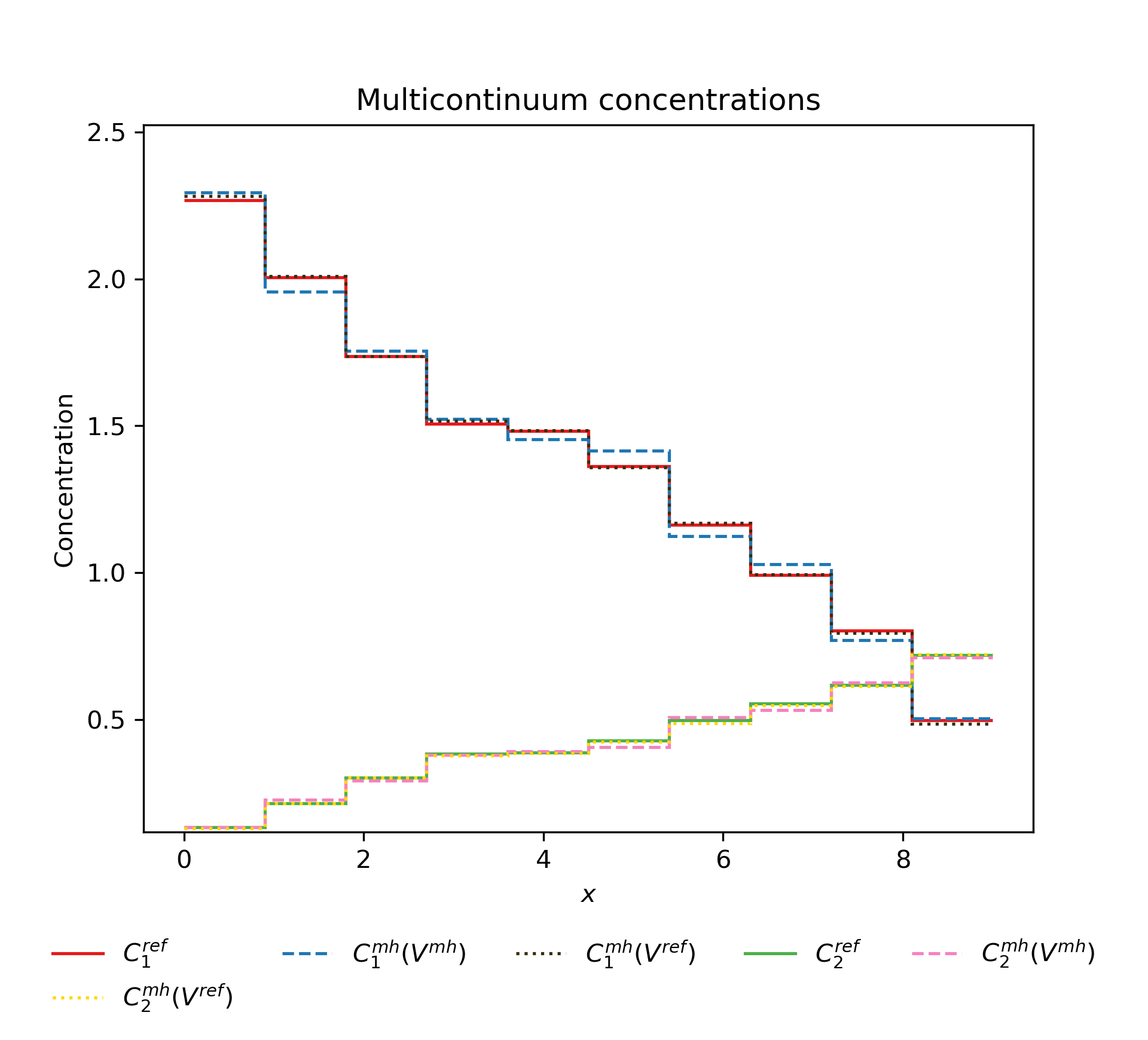}
\caption{Multicontinuum velocity and concentration plots (from left to right) at the final time. Gravity-driven fingering with two continua and heterogeneities.}
\label{fig:coarse_results_problem_4_N_2}
\end{figure}

Figure \ref{fig:coarse_results_problem_4_N_2} presents plots of the multicontinuum velocities and concentrations (from left to right) at the final time. One can observe some changes in both plots caused by the influence of the heterogeneity of the coefficient $\lambda$. We can see that the homogenized velocities follow the trend of the reference multi-continuum velocities, indicating high accuracy. We observe that both homogenized concentrations (based on $V^{\text{ref}}$ and $V^{\text{mh}}$) also approximate the reference concentrations well.

Let us now consider the errors of our multicontinuum approach. For the velocities, we use the relative $L^2$ errors presented in \eqref{eq:errors_U}. In this way, we obtain $e_{V}^{(1)} = 6.880$\% for the first continuum velocity and 
$e_{V}^{(2)} = 7.382$\% for the second continuum velocity at the final time. Therefore, our homogenized velocities can approximate the reference ones with high accuracy.

\begin{table}[hbt!]
\caption{Relative errors of multicontinuum concentrations at the final time for gravity-driven fingering with two continua and heterogeneities.}
\label{tab:C_errors_problem_4_N_2}
\centering
\begin{tabular}{c | c c }
Error $e_{C}^{(i)}$ & $e_{C}^{(1)}$ & $e_{C}^{(2)}$ \\ \hline
$\|C_i^{\text{mh}} (V^{\text{ref}}) - C_i^{\text{ref}} \|_2 / \| C_i^{\text{ref}} \|_2 \times 100\%$ & 0.587\% & 1.020\%\\ 
$\|C_i^{\text{mh}} (V^{\text{mh}}) - C_i^{\text{ref}} \|_2 / \| C_i^{\text{ref}} \|_2 \times 100\%$ & 2.288\% & 2.657\% \\ 
$\|C_i^{\text{mh}} (V^{\text{mh}}) - C_i^{\text{mh}} (V^{\text{ref}}) \|_2  / \| C_i^{\text{mh}} (V^{\text{ref}}) \|_2 \times 100\%$ & 2.334\% &2.656\% \\ 
\end{tabular}
\end{table}

Table \ref{tab:C_errors_problem_4_N_2} presents different types of relative $L^2$ errors for the multicontinuum concentrations at the final time. One can see that all the errors are minor. As expected, the errors of the homogenized concentrations based on the reference velocities $C_i^{\text{mh}}(V^{\text{ref}})$ are smaller than for those based on the homogenized velocities $C_i^{\text{mh}}(V^{\text{mh}})$. The errors between these two homogenized concentrations are also minor and comparable to the errors of $C_i^{\text{mh}}(V^{\text{mh}})$. Therefore, our multicontinuum model can provide high accuracy for heterogeneous dual-continuum media.

\subsubsection{Gravity-driven fingering: three continua with heterogeneities}\label{sec:gravity_fingering_heterogeneous_triple_continuum}


This case represents a gravity-driven fingering problem in a heterogeneous triple-continuum medium. Again, we set the same heterogeneous coefficient $\lambda$ as in the previous problem. The initial condition is taken from the triple-continuum problem with a constant $\lambda$. We set 150 time steps with a step size of $\tau = 4 \cdot 10^{-2}$, the total time is $T = 6$. The CFL number varies between 0.231 and 0.251 over time.

\begin{figure}[hbt!]
\centering
\includegraphics[width=0.49\textwidth]{lmd_problem_4.png}
\includegraphics[width=0.49\textwidth]{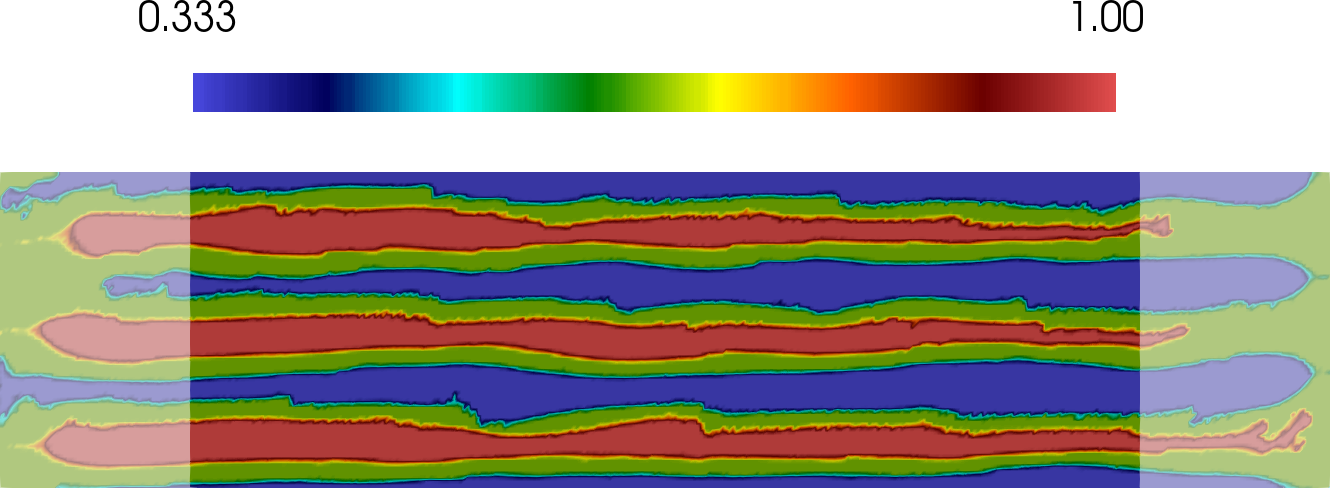}
\caption{Distributions of heterogeneous coefficient $\lambda$ and final fine-scale concentration. Gravity-driven fingering with three continua and heterogeneities.}
\label{fig:c_fine_problem_4_N_3}
\end{figure}

In Figure \ref{fig:c_fine_problem_4_N_3}, we present distributions of the heterogeneous coefficient $\lambda$ and the concentration at the final time (from left to right). We observe the influence of the heterogeneities on the fingers' configuration at the final time. Again, the fingers of the first continuum, corresponding to the densest fluid, are directed to the right in accordance with gravity. In contrast, the third continuum fingers, representing the least dense fluid, are oriented to the left. The fingers of the third continuum experience competing influences from other continua.

\begin{figure}[hbt!]
\centering
\includegraphics[width=0.49\textwidth]{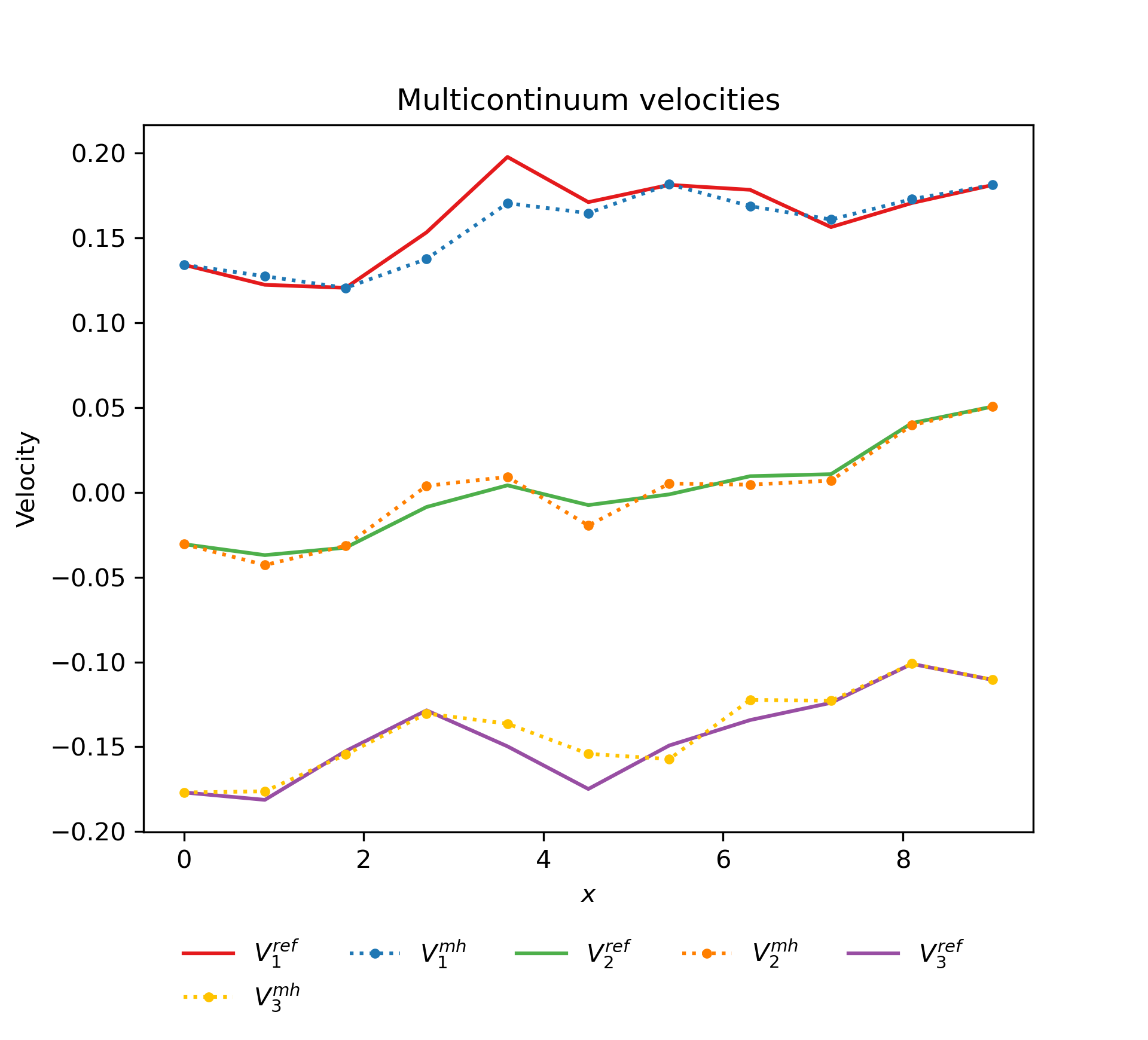}
\includegraphics[width=0.49\textwidth]{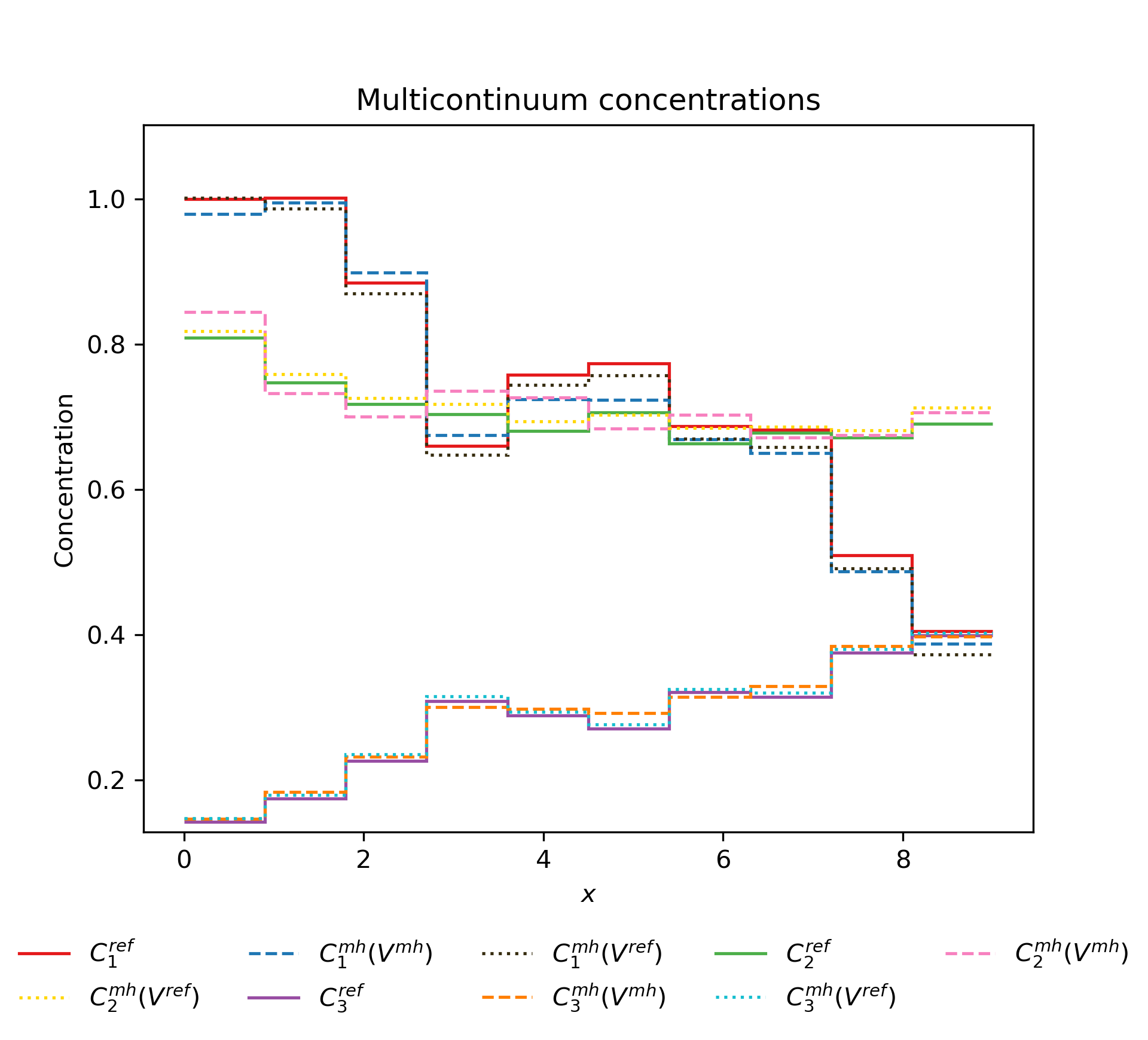}
\caption{Multicontinuum velocity and concentration plots (from left to right) at the final time. Gravity-driven fingering with three continua and heterogeneities.}
\label{fig:coarse_results_problem_4_N_3}
\end{figure}

Figure \ref{fig:coarse_results_problem_4_N_3} presents plots of the multicontinuum velocities and concentrations at the final time (from left to right). We observe the influence of the heterogeneities on the plots. One can observe that the homogenized velocities approximate the reference ones with high accuracy. As in the case with a constant coefficient, we see that the velocity of the second continuum is close to zero due to the symmetry of the initial condition. One can see that both homogenized concentrations ($C_i^{\text{mh}} (V^{\text{ref}})$ and $C_i^{\text{mh}} (V^{\text{mh}})$) accurately follow the trends of the reference multicontinuum concentrations, indicating the high accuracy.

Let us now consider the errors of our multicontinuum solutions. As in the previous triple-continuum problem, we cannot use the relative $L^2$ errors for the multicontinuum velocities, since the second continuum velocity is close to zero on average. Instead, we use the absolute $L^2$ errors and general relative $L^2$ error \eqref{eq:errors_U_N_3}. In this way, we obtain the following absolute errors for each continuum at the final time: $e_{V, \text{abs}}^{(1)} = 0.034$, $e_{V, \text{abs}}^{(2)} = 0.021$, $e_{V, \text{abs}}^{(3)} = 0.029$. For the general relative $L^2$ error, we have $e_{V} = 6.800$\%. Therefore, our homogenized multicontinuum velocities can provide high accuracy in the heterogeneous triple-continuum case.

\begin{table}[hbt!]
\caption{Relative errors of multicontinuum concentrations at the final time for gravity-driven fingering with three continua and heterogeneities.}
\label{tab:C_errors_problem_4_N_3}
\centering
\begin{tabular}{c | c c c }
Error $e_{C}^{(i)}$ & $e_{C}^{(1)}$ & $e_{C}^{(2)}$ & $e_{C}^{(3)}$ \\ \hline
$\|C_i^{\text{mh}} (V^{\text{ref}}) - C_i^{\text{ref}} \|_2 / \| C_i^{\text{ref}} \|_2 \times 100\%$ & 2.394\% & 1.889\% & 1.873\% \\ 
$\|C_i^{\text{mh}} (V^{\text{mh}}) - C_i^{\text{ref}} \|_2 / \| C_i^{\text{ref}} \|_2 \times 100\%$ & 3.394\% & 3.813\% & 3.531\% \\ 
$\|C_i^{\text{mh}} (V^{\text{mh}}) - C_i^{\text{mh}} (V^{\text{ref}}) \|_2  / \| C_i^{\text{mh}} (V^{\text{ref}}) \|_2 \times 100\%$ & 2.655\% & 2.938\% & 2.902\% \\ 
\end{tabular}
\end{table}

Table \ref{tab:C_errors_problem_4_N_3} presents the relative $L^2$ errors for the final-time multicontinuum concentrations at the final time. We see that all the errors are minor. The homogenized concentrations based on the reference velocities $C_i^{\text{mh}} (V^{\text{ref}})$ have smaller errors than those based on the homogenized velocities $C_i^{\text{mh}} (V^{\text{mh}})$, as expected. Also, the errors between $C_i^{\text{mh}} (V^{\text{mh}})$ and $C_i^{\text{mh}} (V^{\text{ref}})$ are minor and smaller than the errors between $C_i^{\text{mh}} (V^{\text{mh}})$ and $C_i^{\text{ref}}$. Therefore, our multicontinuum model can approximate the reference concentrations with high accuracy for the heterogeneous triple-continuum case.

\subsection{Contrast-determined continua}\label{sec:contrast_determined}


In this subsection, we consider the cases where continua are determined by the high contrast of the coefficient $\lambda$, which depends on the fine-scale concentration field. We consider two model problems in dual-continuum media using mixed and Galerkin multicontinuum modeling approaches.

\subsubsection{Viscous fingering}\label{sec:viscous_fingering}


The present case considers a viscous fingering problem in a dual-continuum medium. The fine-scale mathematical model is defined in the right-extended domain $\Omega_{\text{ext}}^{\text{(R)}}$ (see Figure \ref{fig:grids}) and has the following form
\begin{equation}\label{eq:fine_problem_viscous_fingering}
\begin{split}
v= -\lambda(c) \nabla p \quad &\text{in } \Omega_{\text{ext}}^{\text{(R)}},\\
div(v)=0 \quad &\text{in } \Omega_{\text{ext}}^{\text{(R)}},\\
c_t + v \cdot \nabla c = 0 \quad &\text{in } \Omega_{\text{ext}}^{\text{(R)}} \times (0, T],
\end{split}
\end{equation}
where the coefficient $\lambda$ possesses high contrast, i.e., $\frac{\text{max} \lambda}{\text{min} \lambda} \gg 1$, and depends on the fine-scale concentration field $c$ such that
\begin{equation}\label{eq:lambda}
\lambda(c) = \begin{cases}
1000, & c \in I_1,\\
1, & c \in I_2.
\end{cases}
\end{equation}

We complement the system \eqref{eq:fine_problem_viscous_fingering} with the following boundary conditions
\begin{equation}
\begin{gathered}
v \cdot n = g_{\text{in}} \quad \text{on } \Gamma_{\text{L}}^{\text{ext}}, \quad v \cdot n = 0 \quad \text{on } \Gamma_{\text{B}}^{\text{ext}} \cup \Gamma_{\text{T}}^{\text{ext}},\\
p = p_{\text{out}} \quad \text{on } \Gamma_{\text{R}}^{\text{ext}},\\
c = c_{\text{in}} \quad \text{on } \Gamma_{\text{L}}^{\text{ext}},
\end{gathered}
\end{equation} 
where $g_{\text{in}}$ is a normal flux, $p_{\text{out}}$ is an outlet pressure, $c_{\text{in}}$ is an inlet concentration, and $\Gamma_{\text{L}}^{\text{ext}}$, $\Gamma_{\text{R}}^{\text{ext}}$, $\Gamma_{\text{B}}^{\text{ext}}$, and $\Gamma_{\text{T}}^{\text{ext}}$ are the left, right, bottom, and top boundaries of $\Omega_{\text{ext}}^{\text{(R)}}$, respectively, such that $\partial \Omega_{\text{ext}}^{\text{(R)}} = \Gamma_{\text{L}}^{\text{ext}} \cup \Gamma_{\text{R}}^{\text{ext}} \cup \Gamma_{\text{B}}^{\text{ext}} \cup \Gamma_{\text{T}}^{\text{ext}}$.


As a numerical method for solving the flow problem, we again use a mixed finite element method with the lowest-order RT basis functions for the velocity and piecewise constant basis functions for the pressure. Since the velocity in this problem primarily depends on the coefficient $\lambda$ and is less sensitive to slight variations in the concentration $c$ (compared to the previous problems), we use the streamline upwind Petrov-Galerkin (SUPG) method for spatial discretization of the transport problem on a fine grid. The temporal discretization is based on the Backward Euler method. In each time step, we first solve the flow problem and then the transport problem using the computed velocity.


In the multicontinuum approach, we again set $\phi_i^c = \chi_i \psi_i$. Since in the considered case the coefficient $\lambda$ possesses high contrast, we also set $\phi_i^p = \psi_i$. Multicontinuum basis functions for velocity are constructed by solving mixed cell problems. In this way, we have multicontinuum pressure and velocity fields.

We consider two cell problems for multicontinuum velocities. We formulate the first cell problem in a local domain $\omega_l$, consisting of coarse blocks sharing the edge $E_l$. This cell problem considers horizontal fluxes and can be written as follows
\begin{equation}
\begin{split}
\lambda^{-1}(c) \phi_{i}^{v} + \nabla \phi_{i}^{pv} = 0 \quad \text{in } \omega_l,\\
div(\phi_{i}^{v}) = \theta_i^v \psi_i(c) \quad \text{in } \omega_l.
\end{split}
\end{equation}
We complement this cell problem with a boundary condition $\phi_i^{v} \cdot n = 0$ on $\partial \omega_l$ and an additional constraint $\phi_{i}^{v} \cdot n = \psi_{i}(c)$ on $E_l$. We choose $\theta_i^v$ such that it satisfies $\int_{K} \theta_i^v \psi_i(c) = \int_{E_l} \psi_{i}(c)$ for all $K \subset \omega_l$.

The second cell problem is formulated to account for vertical fluxes, i.e., interfaces between continua. Note that the number of these multicontinuum basis functions is equal to the number of the interfaces, i.e., $N - 1 = 1$ in our case. This cell problem is defined in coarse block separately and has the following form
\begin{equation}
\begin{split}
\lambda^{-1}(c) \phi^{w} + \nabla \phi^{pw} = 0 \quad \text{in } K_l,\\
div(\phi^{w}) = \psi_1(c) - \theta^{w, l} \psi_2(c) \quad \text{in } K_l,
\end{split}
\end{equation}
where $\phi^{w} \cdot n = 0$ on $\partial K_l$. Note that $\theta^{w, l}$ is chosen to satisfy $\int_{K_l} \theta^{w, l} \psi_2(c) = \int_{K_l} \psi_1(c)$.


Let us consider the fine-scale simulation. We set some predefined finger configurations as the initial condition $c_0$ for the transport problem. For boundary conditions, we set $g_{\text{in}} = -1$, $c_{\text{in}} = c_0|_{\Gamma_{\text{L}}^{\text{ext}}}$, and $p_{\text{out}} = 0$. To ensure that fingers of both continua are present on all the coarse-grid edges, we run a pre-simulation of 120 time steps with time step $\tau = 1.88 \cdot 10^{-3}$ (chosen so that, by its end, the fingers reach the right boundary of the target domain). We then perform the main simulation for 150 time steps with the same step size. Therefore, the total time is $T = 0.5076$ (including both pre-simulation and main simulation stages). The CFL number varies from 0.181 to 0.420 over time. Note that the coarse-scale simulations will be performed only for the main simulation stage. The pre-simulation stage is used solely to simplify subsequent coarse-scale computations. In general, multicontinuum modeling can be carried out without it.

\begin{figure}[hbt!]
\centering
\includegraphics[width=0.49\textwidth]{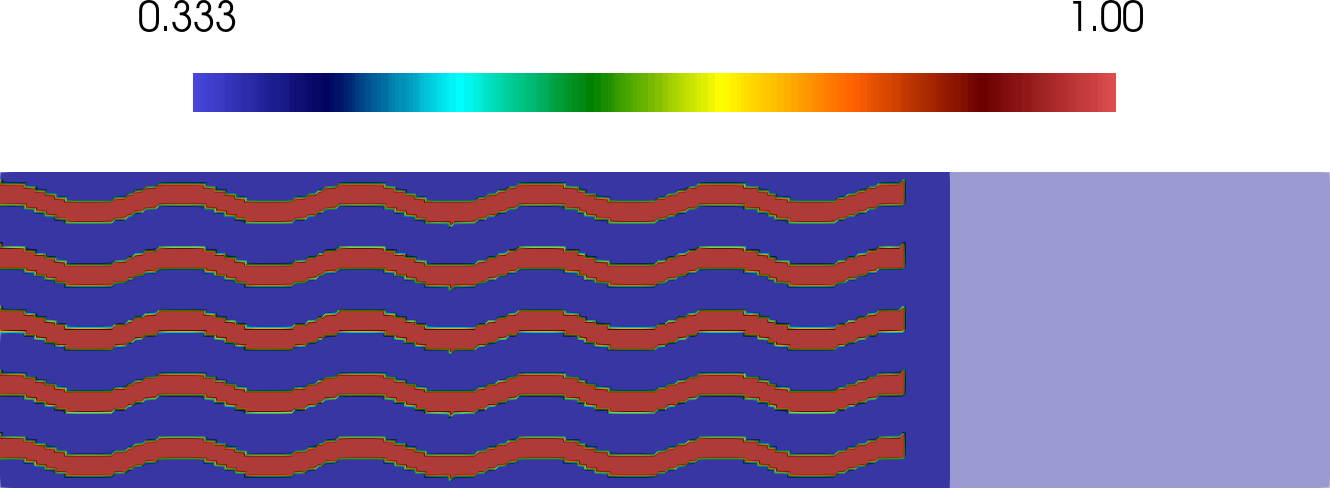}
\includegraphics[width=0.49\textwidth]{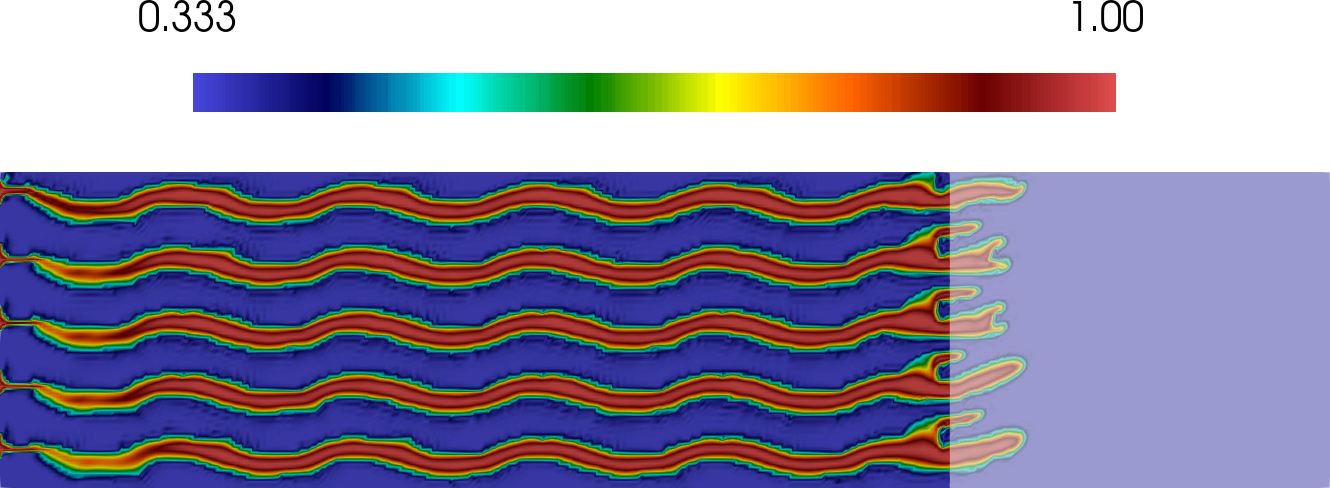}
\caption{Initial and final fine-scale concentration distributions. Viscous fingering.}
\label{fig:c_fine_problem_2_N_2}
\end{figure}

Figure \ref{fig:c_fine_problem_2_N_2} shows distributions of the fine-scale concentration at the initial and final time instances (from left to right). The finger structure evolves significantly during the simulation: the initial fingers deform, and secondary fingers nucleate, consistent with a viscous-fingering instability driven by the coefficient (viscosity) contrast between the two continua.

\begin{figure}[hbt!]
\centering
\includegraphics[width=0.49\textwidth]{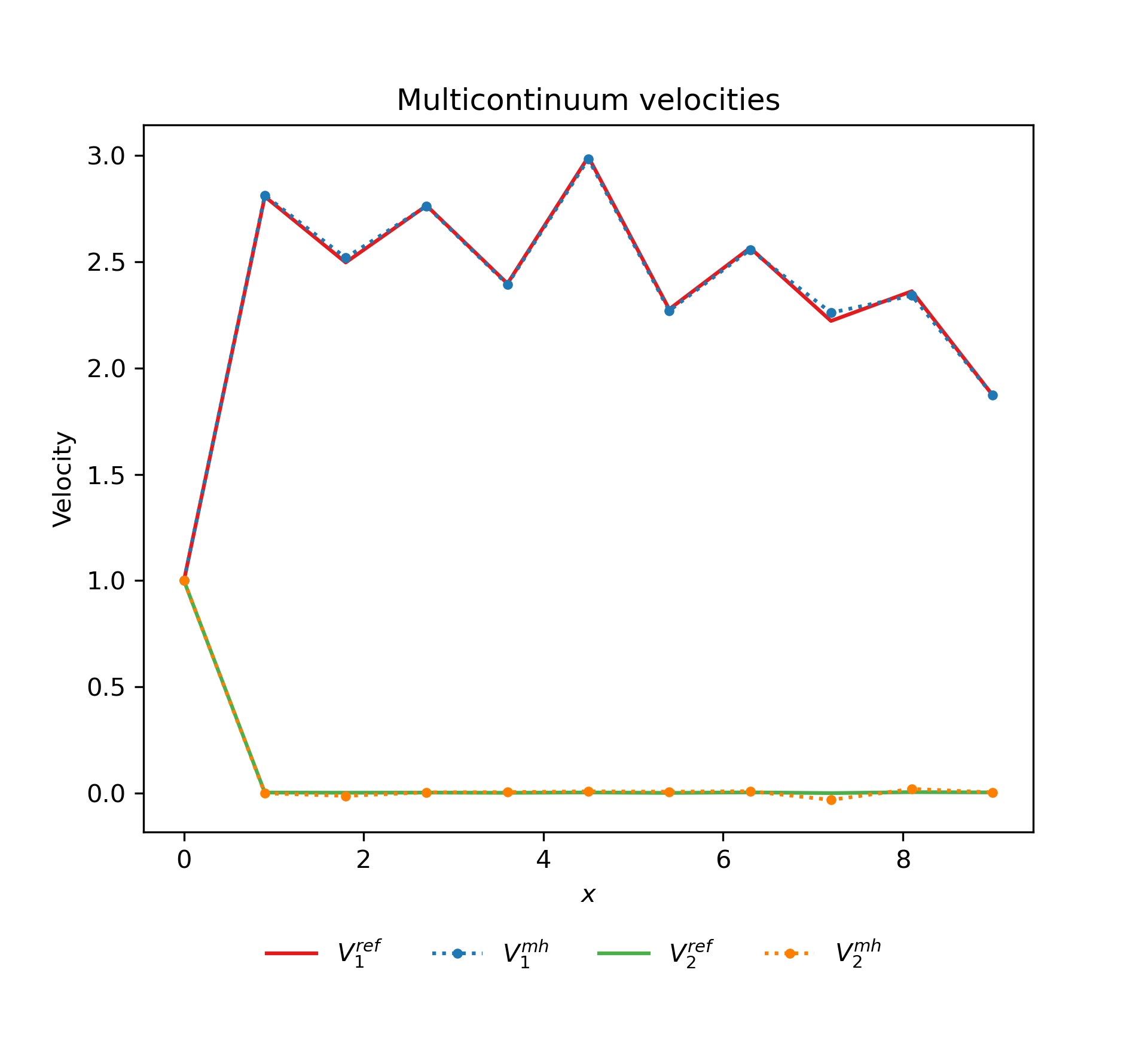}
\includegraphics[width=0.49\textwidth]{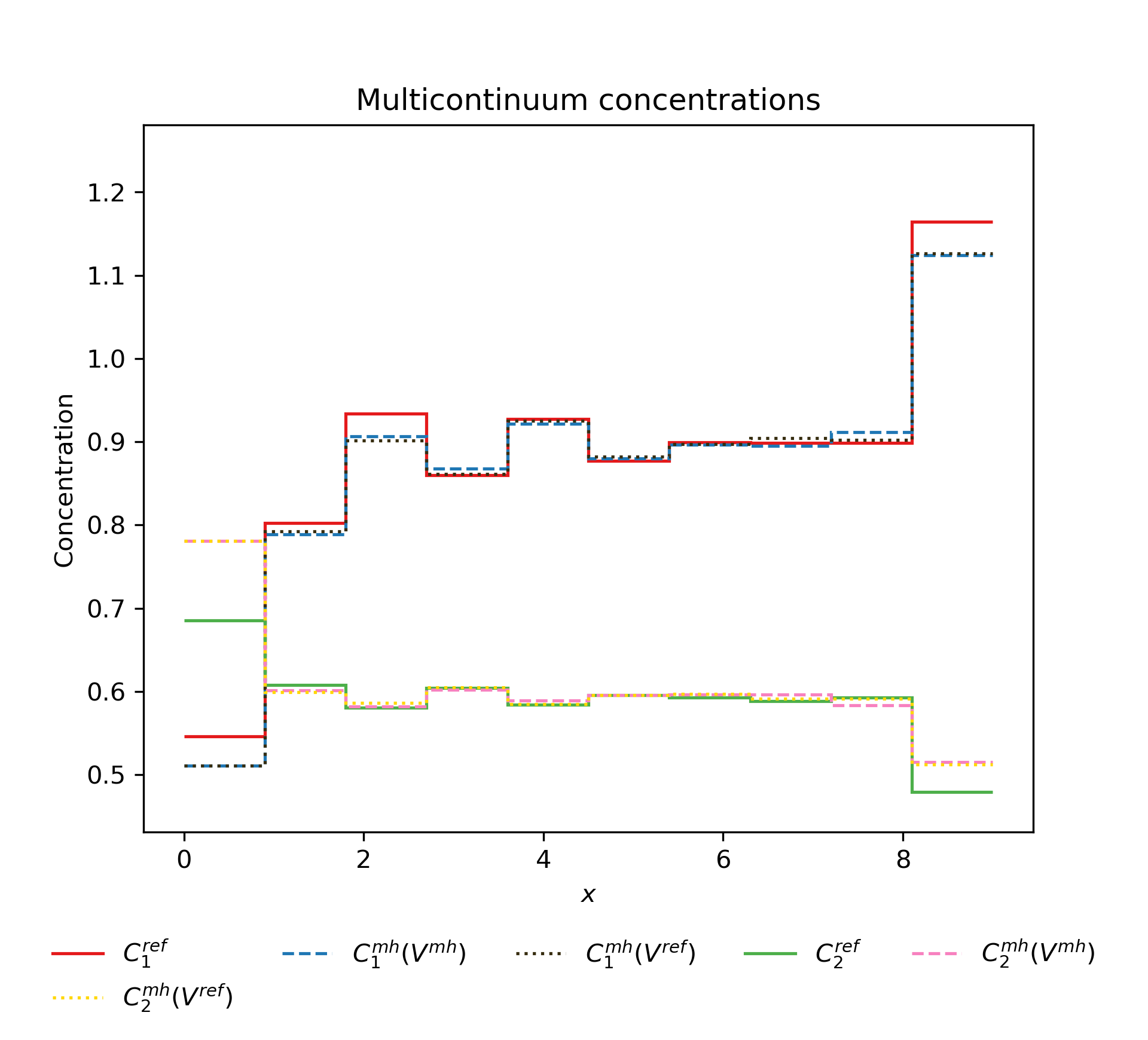}
\caption{Multicontinuum velocity and concentration plots (from left to right) at the final time. Viscous fingering.}
\label{fig:coarse_results_problem_2_N_2}
\end{figure}

In Figure \ref{fig:coarse_results_problem_2_N_2}, we present plots of the multicontinuum velocities and concentrations (from left to right) at the final time. One can see that our homogenized velocities are very similar to the reference multicontinuum velocities, which indicates the high accuracy of the proposed multicontinuum approach. Both homogenized concentrations (based on $V^{\text{ref}}$ and $V^{\text{mh}}$) follow the trends of the reference multicontinuum concentrations, demonstrating the capability to predict concentration dynamics with high accuracy.

Let us consider the errors of our multicontinuum solutions. One can notice that the second continuum velocity is much smaller than the first continuum velocity. However, due to the inflow boundary condition, it is not zero on average. Therefore, we can use the relative $L^2$ errors \eqref{eq:errors_U} for each continuum. We obtain $e_{V}^{(1)} = 0.652$\% for the first continuum and $e_{V}^{(2)} = 3.938$\% for the second continuum at the final time. Thus, our homogenized velocities accurately approximate the reference ones.

\begin{table}[hbt!]
\caption{Relative errors of multicontinuum concentrations at the final time for viscous fingering.}
\label{tab:C_errors_problem_2_N_2}
\centering
\begin{tabular}{c | c c }
Error $e_{C}^{(i)}$ & $e_{C}^{(1)}$ & $e_{C}^{(2)}$ \\ \hline
$\|C_i^{\text{mh}} (V^{\text{ref}}) - C_i^{\text{ref}} \|_2 / \| C_i^{\text{ref}} \|_2 \times 100\%$ & 2.223\%& 5.392\% \\ 
$\|C_i^{\text{mh}} (V^{\text{mh}}) - C_i^{\text{ref}} \|_2 / \| C_i^{\text{ref}} \|_2 \times 100\%$ & 2.264\% & 5.484\% \\ 
$\|C_i^{\text{mh}} (V^{\text{mh}}) - C_i^{\text{mh}} (V^{\text{ref}}) \|_2  / \| C_i^{\text{mh}} (V^{\text{ref}}) \|_2 \times 100\%$ & 0.588\% & 0.620\% \\ 
\end{tabular}
\end{table}

Table \ref{tab:C_errors_problem_2_N_2} presents different types of relative $L^2$ errors of our multicontinuum concentrations at the final time. We see that all the errors are minor. Both homogenized errors $C_i^{\text{mh}} (V^{\text{ref}})$ and $C_i^{\text{mh}} (V^{\text{mh}})$ have similar errors. However, one can see that the errors between $C_i^{\text{mh}} (V^{\text{mh}})$ and $C_i^{\text{mh}} (V^{\text{ref}})$ than between $C_i^{\text{mh}} (V^{\text{mh}})$ and $C_i^{\text{ref}}$. This phenomenon indicates that numerical diffusion plays a significant role in the concentration errors. In general, our multicontinuum model provides accurate approximations of the reference multicontinuum concentrations for the considered viscous fingering problem.

\subsubsection{Interface flattening driven by high-contrast flow}\label{sec:interface_flattening}


In this case, we consider a problem of interface flattening driven by high-contrast flow. Thus, the initial concentration $c_0$ partitions the domain into two parts separated by a wave-like interface. Moreover, the coefficient 
$\lambda(c)$ possesses high contrast \eqref{eq:lambda}, inducing a strongly nonuniform inflow that flattens the interface. The fine-scale mathematical model is defined on the target domain $\Omega$ and can be written as
\begin{equation}
\label{eq:fine_problem_galerkin}
\begin{split}
-div(\lambda(c) \nabla p) =0 \quad &\text{in } \Omega,\\
c_t + v \cdot \nabla c =0 \quad &\text{in } \Omega \times (0, T],
\end{split}
\end{equation}
where $v = -\lambda(c) \nabla p$.

The system \eqref{eq:fine_problem_galerkin} is complemented with the following boundary conditions
\begin{equation}
\begin{gathered}
p = p_{\text{in}} \quad \text{on } \Gamma_{\text{L}}, \quad p = p_{\text{out}} \text{ on } \Gamma_{\text{R}},\\
-\lambda(c) \nabla p \cdot n = 0 \quad \text{on } \Gamma_{\text{B}} \cup \Gamma_{\text{T}},\\
c = c_{\text{in}} \quad \text{on } \Gamma_{\text{L}},
\end{gathered}
\end{equation}
where $p_{\text{in}}$ is an inlet pressure, and $\Gamma_{\text{L}}$, $\Gamma_{\text{R}}$, $\Gamma_{\text{B}}$, and $\Gamma_{\text{T}}$ are the left, right, bottom, and top boundaries of $\Omega$, respectively, such that $\partial \Omega = \Gamma_{\text{L}} \cup \Gamma_{\text{R}} \cup \Gamma_{\text{B}} \cup \Gamma_{\text{T}}$.


In our fine-scale simulations, we use the Galerkin finite element method with standard linear basis functions for the flow problem. For the transport problem, we use the SUPG method for spatial approximation and the Backward Euler method for temporal approximation. Note that, in each time step, we first solve the flow problem and then the transport problem using the computed velocity.


In our multicontinuum modeling approach, we set $\phi_i^c = \chi_i \psi_i$ as the concentration's multicontinuum basis functions. Whereas the pressure basis functions are obtained by solving cell problems in oversampled coarse blocks (RVEs) $K^+$. The oversampled coarse blocks are constructed as layered extensions of coarse blocks $K$.

It should be noted that the cell problems require special treatment to obtain accurate solutions. First, we formulate the cell problems and, accordingly, solve the flow problem on a refined coarse grid. In our simulations, we use the $20 \times 1$ coarse grid for this purpose. The concentrations are computed on the standard $10 \times 1$ coarse grid, as usual. The refined coarse grid for the flow problem is beneficial to obtain accurate multicontinuum velocities defined on the coarse edges $E_l$. We compute the edge-based multicontinuum velocities as the average of the two adjacent block-based velocities (sharing the edge $E_l$).

Second, we allow the oversampled coarse blocks to extend beyond the target domain to the left and right. For the left extension, we periodically extend the concentration field using the initial condition, which is periodic in our example. For the right extension, we reflect the concentration field across the right boundary. In principle, one can use the reflection technique for the left extension as well. It is also possible to solve the cell problems without this extension approach, but in that case, boundary effects may affect the multicontinuum velocities.

We consider both our cell problems in the oversampled coarse blocks constructed with six layers. The first cell problem accounts for gradient effects in each continuum. Considering the problem formulation of the present case, we assume anisotropic flow and ignore the vertical gradient effects. Note that the overflow between continua is accounted for in the second cell problem.
\begin{equation}
\begin{split}
\int_{K^+} \lambda(c) \nabla \bar{\phi}^{p}_i \cdot \nabla q -  \sum_{j,l} \bar{\Gamma}_{ij}^{p, l} \int_{K_l} \psi_j(c)  q  =0,\\
\int_{K_l}  \bar{\phi}^{p}_i \psi_j(c) = \delta_{ij} \int_{K_l} (x_1 - \tilde{x}_{1}) \psi_j(c),
\end{split}
\end{equation}
where $\tilde{x}_1$ is defined by $\int_{K_{l_0}} (x_1 - \tilde{x}_{1}) = 0$.

The second cell problem accounts for different averages in each continuum and can be written as follows
\begin{equation}
\begin{split}
\int_{K^+} \lambda(c) \nabla \phi_i^p\cdot \nabla q -
\sum_{j,l} \Gamma_{ij}^{p,l} \int_{K_l} \psi_j(c)  q =0,\\
\int_{K_l}  \phi_i^p \psi_j(c) = \delta_{ij} \int_{K_l} \psi_j(c).
\end{split}
\end{equation}

\begin{figure}[hbt!]
\centering
\includegraphics[width=0.49\textwidth]{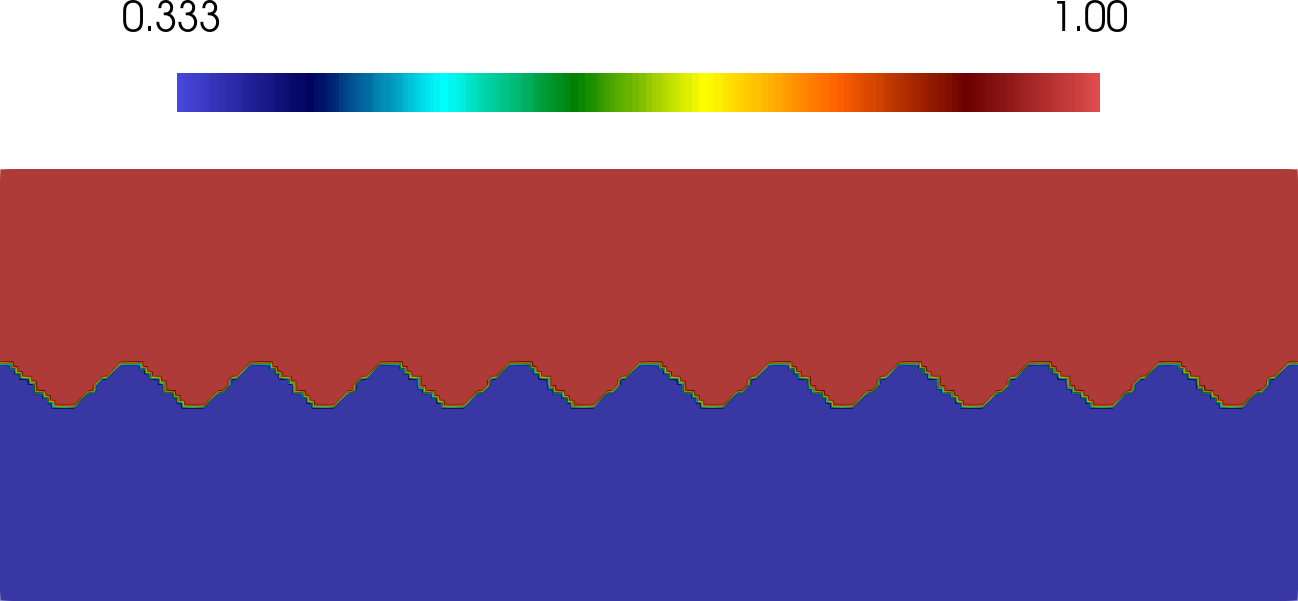}
\includegraphics[width=0.49\textwidth]{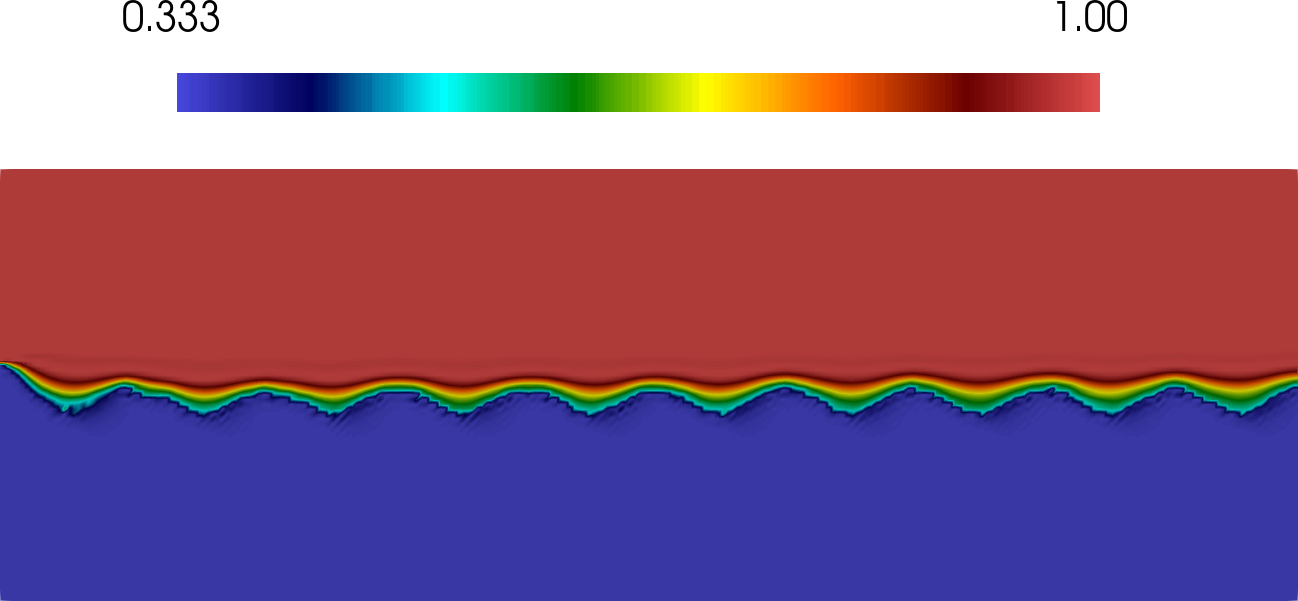}
\caption{Initial and final fine-scale concentration distributions. Interface flattening driven by high-contrast flow.}
\label{fig:c_fine_problem_1_N_2}
\end{figure}

Let us now consider the fine-scale results. We simulate for 1000 time steps with step size $\tau = 1.5 \cdot 10^{-4}$, the total time is $T = 0.15$. For the boundary conditions, we set $p_{\text{in}} = 1$, $p_{\text{out}} = 0$, and $c_{\text{in}} = c_{0}|_{\Gamma_L}$. The CFL number varies from 0.745 to 0.858 over time. Figure \ref{fig:c_fine_problem_1_N_2} shows initial and final distributions of the concentration field (from left to right). One can observe that the wave-like interface has flattened out due to the influence of the high-contrast flow directed from left to right.

\begin{figure}[hbt!]
\centering
\includegraphics[width=0.49\textwidth]{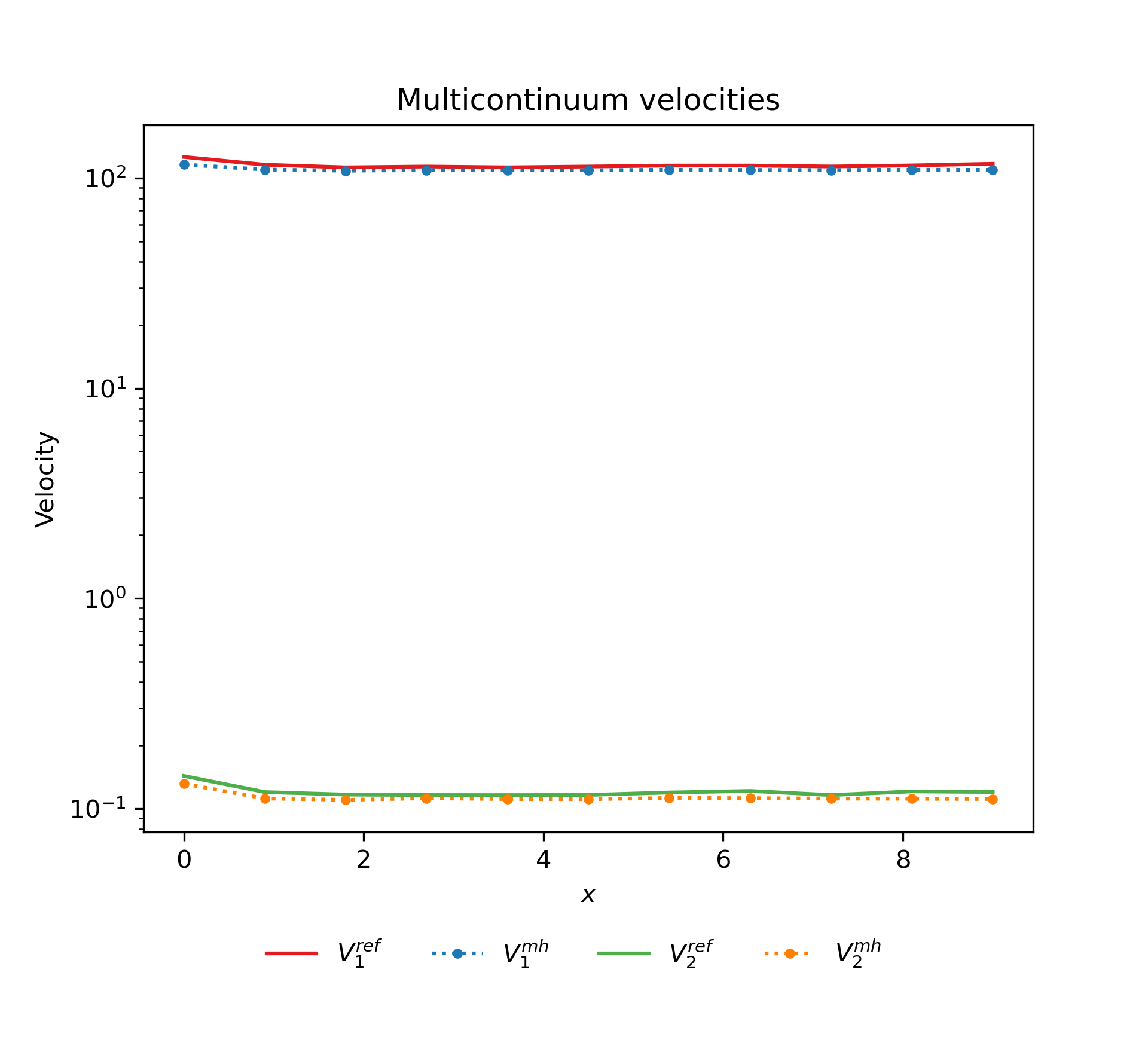}
\includegraphics[width=0.49\textwidth]{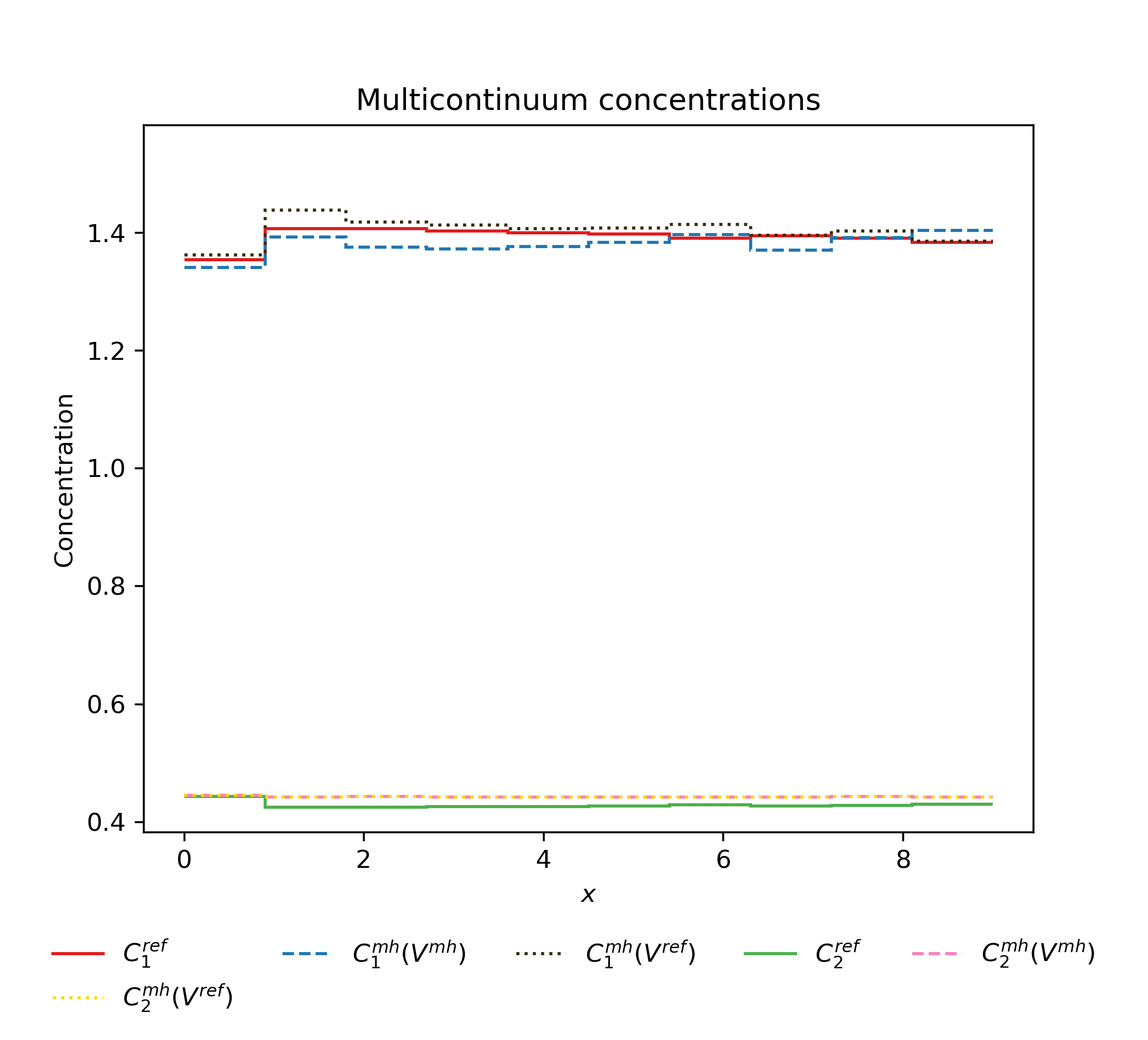}
\caption{Multicontinuum velocity and concentration plots (from left to right) at the final time. Interface flattening driven by high-contrast flow.}
\label{fig:coarse_results_problem_1_N_2}
\end{figure}

In contrast to the previous cases, we use a different number of time steps for coarse-scale modeling. We simulate for 100 time steps with step size $\tau_{\text{coarse}} = 1.5 \cdot 10^{-3}$. In Figure \ref{fig:coarse_results_problem_1_N_2}, we present plots of the multicontinuum velocities and concentrations (from left to right) at the final time. We see that the velocities of the first and second continua differ significantly. Our homogenized velocities closely approximate the reference multicontinuum velocities. Likewise, the homogenized concentrations $C_i^{\text{mh}}(V^{\text{mh}})$ and $C_i^{\text{mh}}(V^{\text{ref}})$ exhibit high accuracy in approximating the reference concentrations.

Let us now consider the errors of our multicontinuum solutions. We compute the relative $L^2$ errors \eqref{eq:errors_U} for the final-time multicontinuum velocities. In this way, we have $e_{V}^{(1)} = 4.997$\% for the first continuum and $e_{V}^{(2)} = 6.294$ \% for the second continuum. Therefore, our multicontinuum model can provide high accuracy in approximating the reference multicontinuum velocities. 

\begin{table}[hbt!]
\caption{Relative errors of multicontinuum concentrations at the final time for interface flattening driven by high-contrast flow.}
\label{tab:C_errors_problem_1_N_2}
\centering
\begin{tabular}{c | c c }
Error $e_{C}^{(i)}$ & $e_{C}^{(1)}$ & $e_{C}^{(2)}$ \\ \hline
$\|C_i^{\text{mh}} (V^{\text{ref}}) - C_i^{\text{ref}} \|_2 / \| C_i^{\text{ref}} \|_2 \times 100\%$ & 1.052\% & 3.477\%\\ 
$\|C_i^{\text{mh}} (V^{\text{mh}}) - C_i^{\text{ref}} \|_2 / \| C_i^{\text{ref}} \|_2 \times 100\%$ & 1.451\% & 3.483\% \\ 
$\|C_i^{\text{mh}} (V^{\text{mh}}) - C_i^{\text{mh}} (V^{\text{ref}}) \|_2  / \| C_i^{\text{mh}} (V^{\text{ref}}) \|_2 \times 100\%$ & 2.135\% &0.036\% \\ 
\end{tabular}
\end{table}

Table \ref{tab:C_errors_problem_1_N_2} presents relative $L^2$ errors for the multicontinuum concentrations at the final time. We observe that the errors for $C_i^{\text{mh}}(V^{\text{ref}})$ are smaller than those for $C_i^{\text{mh}}(V^{\text{mh}})$, as the former reflects our irreducible error. Moreover, the errors between $C_i^{\text{mh}}(V^{\text{mh}})$ and $C_i^{\text{mh}}(V^{\text{ref}})$ are smaller than those between $C_i^{\text{mh}}(V^{\text{mh}})$ and $C_i^{\text{ref}}$, since the former is less influenced by numerical diffusion. Overall, the results show that our multicontinuum model accurately predicts the dynamics of the multicontinuum concentrations.

\section{Conclusions}\label{sec:conclusions}

In this paper, we have proposed a multicontinuum homogenization approach for nonlinear problems involving intrinsically dynamically evolving multiscale media, where the intrinsic variable for multiscale media evolution is given by a transport equation. In contrast to previous approaches, where continua are stationary and defined by the contrast of properties, the proposed approach treats continua as dynamic and defined by one of the fine-scale functions. Based on this concept, we have developed a general formalism for macroscopic modeling of complex nonlinear problems involving dynamically evolving media. As an example, we have considered a fingering problem, where fine-scale concentration defines continua. We have presented the derivation of the macroscopic models based on Galerkin and mixed multicontinuum approaches. In both approaches, we first construct multicontinuum expansions over macroscopic variables and then formulate cell problems to determine the multicontinuum basis functions. After that, we derive the corresponding macroscopic equations. In the numerical results, we have considered cases of gravity-determined and contrast-determined continua. In the former, we have addressed gravity-driven fingering problems in dual- and triple-continuum media, including the cases with heterogeneous coefficients. In the latter, we studied the problems of viscous fingering and interface flattening driven by high-contrast flow in dual-continuum media. Depending on the problem, we employed mixed and Galerkin multicontinuum approaches. The numerical results demonstrate that the macroscopic models, constructed by the proposed multicontinuum approach, accurately represent coarse-scale solution fields for complex nonlinear problems.

\appendix\section{Another derivation}\label{sec:appendix}

In this appendix, we present another derivation of the multicontinuum concentration equations
\begin{align*}
\partial_{t}X_{s}(t,x) & =\vec{v}(X_{s}(t,x))\\
X_{s}(s,x) & =x.
\end{align*}
We consider a fixed rectangular subdomain $[x_{j},x_{j+1}]\times Y$ (see Figure \ref{fig:grids}), for simplicity. The results can be generalized for general macroscopic coarsening.
For $c_{t}+v\cdot\nabla c=0$ in $[x_{j},x_{j+1}]$
\[
\partial_{t}c(t,X_{0}(t,x))=c_{t}(t,X_{0}(t,x))+v(X_{0}(t,x))\cdot\nabla c(t,X_{0}(t,x))=0.
\]

We define $\psi_{i}(t,x)$ be functions such that $\psi_{i}(0,x)=I_{K_{i}}(x)$
and
$\psi_{i}(t,X_{0}(t,x))=\psi_{i}(0,X_{0}(0,x))=I_{K_{i}}(x)$. We also
define 
\[
\Omega_{i}(t)=\{x|\psi_{i}(t,x)=1\}.
\]
 Next, we have
\begin{align*}
\partial_{t}\int_{\Omega_{i}(t)\cap([x_{j},x_{j+1}]\times Y)}cdx & =\int_{\Omega_{i}(t)\cap([x_{j},x_{j+1}]\times Y)}\partial_{t}c+\nabla\cdot(\tilde{v}c)dx\\
 & =\int_{\Omega_{i}(t)\cap([x_{j},x_{j+1}]\times Y)}\partial_{t}c+(\nabla\cdot (v c))dx+\int_{\Omega_{i}(t)\cap([x_{j},x_{j+1}]\times Y)}\nabla\cdot((\tilde{v}-v)c) dx\\
 & =\int_{\Omega_{i}(t)\cap(\{x_{j},x_{j+1}\}\times Y)}c(\tilde{v}-v)\cdot n_{i}dx=-\int_{\Omega_{i}(t)\cap(\{x_{j},x_{j+1}\}\times Y)}v c\cdot n_{i}dx,
\end{align*}
where 
\[
\tilde{\ensuremath{v}}(z)=\begin{cases}
v(z) & z_{1}\notin\{x_{j},x_{j+1}\}\\
(0,v_{2}(z)) & z_{1}\in\{x_{j},x_{j+1}\}
\end{cases}
\]
Therefore, we have 
\[
\partial_{t}\int_{\Omega_i(t)\cap([x_{j},x_{j+1}]\times Y)}cdx=-\int_{\Omega_{i}(t)\cap(\{x_{j},x_{j+1}\}\times Y)}c(v)\cdot n_{i}dx.
\]

\bibliographystyle{unsrt}
\bibliography{lit}

\end{document}